\documentclass[authoryear,review,onefignum,onetabnum]{siamart190516}

\usepackage[semicolon]{natbib}
\allowdisplaybreaks
\usepackage{adjustbox}
\usepackage{enumerate}
\usepackage{bbm}
\usepackage{float}
\usepackage[shortlabels]{enumitem}

\usepackage{bm}
\usepackage{booktabs}
\usepackage{multirow}
\usepackage{threeparttable}
\newcommand\Item[1][i]{%
	\ifx\relax#1\relax  \item \else \item[#1] \fi
	\abovedisplayskip=0pt\abovedisplayshortskip=0pt~\vspace*{-\baselineskip}}

\usepackage{lipsum}
\usepackage{amsfonts}
\usepackage{graphicx}
\usepackage{epstopdf}
\usepackage{algorithmic}
\ifpdf
  \DeclareGraphicsExtensions{.eps,.pdf,.png,.jpg}
\else
  \DeclareGraphicsExtensions{.eps}
\fi

\newsiamremark{remark}{Remark}
\newsiamremark{hypothesis}{Hypothesis}
\crefname{hypothesis}{Hypothesis}{Hypotheses}
\newsiamthm{claim}{Claim}

\headers{Well-posedness and Stability Analysis of Generalized SVM}{N. Ning and J. Wu}

\title{Well-posedness and Stability Analysis of Two Classes of Generalized Stochastic Volatility Models}

\author{Ning Ning\thanks{Department of Statistics, University of Michigan, Ann Arbor, MI 48109
  (\email{patning@umich.edu}).}
\and 
Jing Wu\thanks{School of Mathematics, Sun Yat-sen University, Guangzhou, Guangdong 510275, People's Republic of China
  (\email{wjjosie@hotmail.com}).}}

\usepackage{amsopn}

\def\bt{\begin{theorem}}
\def\et{\end{theorem}}
\def\bl{\begin{lemma}}
\def\el{\end{lemma}}
\def\br{\begin{remark}}
\def\er{\end{remark}}
\def\be{\begin{equation}}
\def\ee{\end{equation}}
\def\ce{\begin{equation*}}
\def\de{\end{equation*}}

\def\cC{{\mathcal C}}

\def\mR{{\mathbb R}}

\def\<{{\langle}}
\def\>{{\rangle}}

\def\e{{\epsilon}}

\newtheorem*{Main Theorem}[theorem]{Main Theorem}{\normalfont\bfseries}{\itshape}
{\bfseries}{\itshape} 
{\bfseries}{\itshape} 
\newtheorem{condition}{Condition}[section]{\bfseries}{\itshape} 
\ifpdf
\hypersetup{
  pdftitle={An Example Article},
  pdfauthor={D. Doe, P. T. Frank, and J. E. Smith}
}
\fi

\nolinenumbers
\begin{document}

\maketitle

\begin{abstract}
In this paper, to cope with the shortage of sufficient theoretical support resulted from the fast-growing quantitative financial modeling, we investigate two classes of generalized stochastic volatility models, establish their well-posedness of strong solutions, and conduct the stability analysis with respect to small perturbations. In the first class, a multidimensional path-dependent  process is driven by another multidimensional path-dependent process. The second class is a generalized one-dimensional stochastic volatility model with H\"older continuous coefficients. What greatly differentiates those two classes of models is that both the process and its correlated driving process have their own subdifferential operators, whose one special case is the general reflection operators for multi-sided barriers. Hence, the models investigated fully cover various newly explored variants of stochastic volatility models whose well-posedness is unknown, and naturally serve as the rigorous mathematical foundation for new stochastic volatility model development in terms of multi-dimension, path-dependence, and multi-sided barrier reflection.  
\end{abstract}



\section{Introduction}
\label{sec1a}\vspace{-4mm}
Stochastic volatility is one of the main concepts widely used in mathematical finance to deal with the endemic time-varying volatility and co-dependence found in financial markets.
Stochastic volatility models since its invention have been widely used to evaluate derivative securities such as options, with the characteristic that the variance of a stochastic process is itself randomly distributed. Various extensions of stochastic volatility models for different purposes have been proposed in recent years with the fast-growing quantitative financial modeling in the past decade. 
However, a shortage of sufficient theoretical support in terms of the existence and uniqueness of a (strong) solution of the proposed models comes along. 

To cope with that, in this paper, we consider two classes of generalized stochastic volatility models, establish their well-posedness, and conduct stability analysis. The first class is the multi-dimensional path-dependent system \eqref{eqn:OriginalsystemofFBSDEs}, where a $d_2$-dimensional path-dependent $Y$ process is driven by a $d_1$-dimensional path-dependent $X$ process. The second class is a generalized one-dimensional stochastic volatility model with H\"older continuous coefficients \eqref{eqn:holder}. What greatly differentiates those two classes of models is that both the $X$ and $Y$ processes have their own subdifferential operators, whose one special case is the general reflection operators for multi-sided barriers, because of which the models under investigation are called stochastic variational inequalities (SVI).   

For illustrative purpose, we consider a simplified one-dimensional path-dependent version of  \eqref{eqn:OriginalsystemofFBSDEs} without control as follows
\begin{equation}\label{sec:OriginalsystemofFBSDEs1d}
\left\{\begin{array}{lll}
 X_t\in x_{0}+\int_{0}^tb(s,X(s))ds+\int_{0}^t\sigma(s,X(s))d\widehat{W}_s-\int_{0}^t\partial\psi_1(X_s)ds, \\
       \\
Y_t\in y_{0}+\int_{0}^t \alpha(s,X(s),Y(s))ds+\int_{0}^t\beta(s,X(s),Y(s))dB_s-\int_{0}^t\partial\psi_2(Y_s)ds,
 \end{array}\right.
\end{equation}
where the path $X(t):=X_{t\wedge\cdot}$ up to time $t$, $\widehat{W}:=\sqrt{1-\rho^2}W+\rho B$ for $W$ and $B$ being two independent one-dimensional Brownian motions with $d\langle\widehat{W},B\rangle_t=\rho dt$ for $|\rho|\leq 1$. Apparently, \eqref{sec:OriginalsystemofFBSDEs1d} covers all the classical types of 
stochastic volatility models and path-dependent models, and it also covers the Heston-type stochastic path-dependent volatility  model proposed in \cite{cozma2018strong} (as well as local maximum stochastic volatility model proposed in \cite{bain2019calibration}) whose well-posedness is unknown,
\begin{align*}
dS_t&=\mu(t,S_t,M_t)S_tdt+\sqrt{V_t}\sigma(t,S_t,M_t)S_tdW_t,\\
dV_t&=\kappa(\theta-V_t)dt+\xi\sqrt{X_t}dW_t^{V},
\end{align*}
where $\sigma$ is a local volatility function depending on the running maximum $M_t:=\sup_{0\leq u\leq t}S_u$, and  $d\langle W ,W^{V}\rangle_t=\rho dt$ for $|\rho|\leq 1$.

%

Reflection factors on stochastic differential equations (SDEs) have wide application and a long history in financial mathematics with great contributions from the pioneer works of N. El Karoui since $1970$s, see \cite{elkaroui1975processus}. 
For economic
dynamics, reflected SDEs was used for the target zone models of the currency
exchange rate (see, for example, \cite{krugman1991target, bertola1993stochastic}). In a regulated financial market,
government regulations lead the spot foreign exchange
(FX) rate processes, the domestic interest rate processes,
and the goods or services (for
instance, grain, water, gas, electricity supply and other
important materials or services for a country), because of which reflected SDEs can be applied realistically and appropriately (see, for example, \cite{bo2011conditional,bo2011some,bo2013conditional}). \eqref{sec:OriginalsystemofFBSDEs1d}
not only extends all the classical reflected SDEs to handle multi-sided barriers, but also covers new models such as reflected stochastic local volatility model in its generalized skew stochastic local volatility model proposed in \cite{ding2020markov} (as well as the reflected stochastic volatility model proposed therein) whose well-posedness is unknown, by taking the special form $\psi_1(X_t)=(2p-1)\mathbbm{1}_{\{X_t\geq a\}}$,
\begin{align*}
dS_t&=\gamma(S_t,X_t)dt+m(X_t)\gamma(S_t)dW_t^{(1)},\\
dX_t&=\mu(X_t)dt+\sigma(X_t)dW_t^{(2)}+(2p-1)dL_t^X(a),
\end{align*}
where $d\langle W^{(1)},W^{(2)}\rangle_t=\rho dt$ for $|\rho|\leq 1$, and $L_t^X(a)$ is the symmetric local time of $X$ at the point $a$, and $p=0$ or $1$ for the $X$ process being the reflected diffusion at the value $a$.

Following the new trend in financial mathematics, a control process belonging to the set of predictable processes and taking values in a compact separable metric space, is embed in both the drift function and the diffusion function of the $Y$ process of both two classes of models under investigation. This control process equips the proposed models the applicability in stochastic control problems, such as the super-replicate valuation problem using the uncertain volatility models with stochastic bounds in \cite{fouque2018uncertain}. 
We further followed \cite{fouque2018uncertain} in conducting the stability analysis of the SVI  systems \eqref{eqn:OriginalsystemofFBSDEs} and \eqref{eqn:holder} by perturbing the systems with a small positive parameter $\epsilon$. Asymptotic analyses were conducted on the perturbed systems to explore their limiting behaviors as $\epsilon$ goes to zero.  In financial mathematics, stochastic volatility models with a small parameter is a typical setup (see, for example, \cite{fouque2000derivatives,fouque2011multiscale}), which may function on the driving volatility process ($X$ process in the current setting) resulting in slow-moving effects.

 Well-posedness for the two classes of models has to be established by different methods due to very different model setups. On proving the well-posedness of the multidimensional SVI system \eqref{eqn:OriginalsystemofFBSDEs}, we used the method of Euler scheme for any duration $T$. To handle the path-dependent effects, we extensively applied the functional It\^o formula that was introduced by \cite{dupire2019functional}.
When it comes to the one-dimensional model with H\"older continuous coefficients \eqref{eqn:holder}, we established its wellposedness by means of the Moreau-Yosida regularization approximation method which was used in \cite{asiminoaei1997approximation} with Lipschitz continuous coefficients. Analogous techniques can be used in handling other problems, see for example, \cite{ren2016approximate} on approximating continuity and the support of reflected SDEs, \cite{ren2013optimal} on reflected SDEs with jumps and its associated optimal control problems, \cite{wu2018limit} on limit theorems and the support of SDEs with oblique reflections on nonsmooth domains.

The rest of the paper is organized as follows. 
In Section \ref{sec:general_model}, we analyze the multi-dimensional path-dependent SVI system \eqref{eqn:OriginalsystemofFBSDEs}, where the well-posedness of the  $X$ and $Y$ processes is established in Section \ref{sec:general_model_wellposedness_X}
and Section \ref{sec:general_model_wellposedness_Y} respectively. Next we considered a perturbed version of \eqref{eqn:OriginalsystemofFBSDEs} with a small positive parameter $\epsilon$, and showed that the perturbed $X^\epsilon$ and $Y^\epsilon$ processes converge to the $X$ and $Y$  processes in Sections \ref{sec:general_model_asymptotic_X} and \ref{sec:general_model_asymptotic_Y} respectively. 
In Section \ref{sec:holder}, we investigate the one-dimensional model with H\"older continuous coefficients \eqref{eqn:holder}, whose well-posedness is established in Section \ref{sec:holder_wellposedness} and whose stability analysis is conducted in Section \ref{sec:holder_asymptotic}. In the sequel, $C$ stands for a constant which may change line by line.

\section{Multi-dimensional Path-dependent SVI  System}\label{sec:general_model}

In this section, our investigation is based on the following general multi-dimensional path-dependent system of  stochastic variational inequalities (SVI):
\begin{equation}\label{eqn:OriginalsystemofFBSDEs}
\left\{\begin{array}{lll}
 X_t\in & x_{0}+\int_{0}^tb(s,X(s))ds+\int_{0}^t\sigma_1(s,X(s))dW_s+\int_{0}^t\sigma_2(s, X(s)) dB_s\\
 &-\int_{0}^t\partial\psi_1(X_s)ds, \\
       \\
Y_t\in & y_{0}+\int_{0}^t \alpha(s,X(s),Y(s), q_s)ds+\int_{0}^t\beta(s,X(s),Y(s), q_s)dB_s\\
&-\int_{0}^t\partial\psi_2(Y_s)ds.
 \end{array}\right.
\end{equation}
Here, $X_t\in \mathbb R^{d_1}$ denotes the status of $X$ at time $t\in[0,T]$; $b$, $\sigma_1$, and $\sigma_2$ are measurable functions on $\mathbb R^+\times\cC(\mathbb R^+;\mathbb R^{d_1})$  depending on the path $X(t):=X_{t\wedge\cdot}$ up to time $t$, valued in $\mathbb R^{d_1}$, $\mathbb R^{{d_1}\times d_W}$, and $\mathbb R^{{d_1}\times d_B}$, respectively; $W$ and $B$ are two independent $d_W$-dimensional and $d_B$-dimensional standard Brownian motions on a complete filtered probability space $(\Omega,\mathcal F,\{\mathcal F_t;t\geq0\},\mathbb P)$. 
We call $\nu:=(\Omega,\mathcal F,\{\mathcal F_t;t\geq0\},\mathbb P, W, B)$ a reference system, based on which, denote $\mathcal{A}_{\nu}$ as the set of admissible controls that is the set of $(\mathcal F_t)$-predictable and $\mathbb U$-valued processes. $Y_t\in \mathbb R^{d_2}$ denotes the status of $Y$ at time $t\in[0,T]$; $q$ is the control process belonging to the set of predictable processes and taking values in a compact separable metric space $\mathbb U$; $\alpha$ and $\beta$ are measurable functions on $\mathbb R^+\times\cC(\mathbb R^+;\mathbb R^{d_1})\times\cC(\mathbb R^+;\mathbb R^{d_2}) \times\mathbb U$, valued in $\mathbb R^{d_2}$ and $\mathbb R^{{d_2}\times d_B}$ respectively, depending on both paths $X(t)$ and $Y(t)$ as well as the control process $q$. 

For $i=1,2$, $\psi_i$ is a proper, convex, and lower-semicontinuous function on $\mathbb R^{d_i}$, with its effective domain
$$
D_i:=\{x\in\mathbb R^{d_i}: \psi_i(x)<\infty\},
$$
and its subdifferential operator 
$$\partial\psi_i(x):=\{z\in\mathbb R^{d_i}; \langle x'-x,z\rangle\leq\psi_i(x')-\psi_i(x), \forall x'\in\mathbb R^{d_i}\},$$
where $\langle\cdot,\cdot\rangle$ denotes the inner product.
Theories on subdifferential operators (see, \cite{rockafellar1970maximal}) indicate that $\partial\psi_i(x)$ is closed and convex for every $x\in\mathbb R^{d_i}$, satisfying that 
$$
\langle x-x',z-z'\rangle\geq0
$$
for any $x, ~x'\in\mR^{d_i}$, $z\in\partial\psi_i(x)$, and $z'\in\partial\psi_i(x')$; 
 $\partial\psi_i$ is maximal monotone, that is, if $x, ~z\in\mR^{d_i}$ satisfying
$$\langle x-x',z-z'\rangle\geq0
$$
for any $x'\in\mathbb R^{d_i}$ and $z'\in\partial\psi_i(x')$, then $z\in\partial\psi_i(x)$.

\begin{condition}\label{X}
For the $X$ process in the SVI  system \eqref{eqn:OriginalsystemofFBSDEs}, we impose the following conditions:
\begin{itemize}
\item  $b(t, x)$ and $\sigma_i(t, x)$ are continuous in $t$, and satisfy
\begin{align}
\label{eqn:modelX}
&\langle b(t, x(t))-b(t, x'(t)),x_t-x'_t\rangle\leq 0,\quad &\forall x, x'\in\cC(\mathbb R^+;\mathbb R^{d_1}),\nonumber\\
&|b(t, x(t))-b(t, x'(t))|\leq l_0(t) \|x-x'\|_t^{\frac12+\alpha}, \quad & \text{for some } \alpha\in[0,1/2],\\
&\|\sigma_i(t, x(t))-\sigma_i(t, x'(t))\|\leq l_i(t)\|x-x'\|_t, \quad & i=1, 2,\nonumber
\end{align}
where $l_i(\cdot)\in L^2([0,T])$ for $i=0,1,2$ and $\|z\|_t:=\sup_{s\leq t}|z_s|$.
\item $0\in\mathrm{Int}(D_1)$ and $\psi_1\geq\psi_1(0)\equiv0$.
\end{itemize}
\end{condition}

\begin{condition}
\label{Y}
For the $Y$ process in the SVI  system \eqref{eqn:OriginalsystemofFBSDEs}, we impose the following conditions:
\begin{itemize}
\item $\lambda_1\leq q_t\leq \lambda_2$.
\item For $\|x\|_t\leq R$ and $L_R(t)$ being locally square integrable,
\begin{align*}
\left|\alpha(t, x(t), y(t),q_t)-\alpha(t, x(t), y'(t),q_t)\right|\leq & L_R(t) \|y-y'\|_t,\\
\left\|\beta(t, x(t), y(t),q_t)-\beta(t, x(t), y'(t),q_t)\right\|\leq & L_R(t) \|y-y'\|_t.
\end{align*}
\item  $\alpha(\cdot,\cdot,\eta,\cdot)$ and $\beta(\cdot,\cdot,\eta,\cdot)$ are continuous in $\mR^+\times\cC(\mR^+;\mR^{d_1})\times\mathbb{U}$, for $\eta\in\cC(\mR^+;\mR^{d_2})$.
\item $0\in\mathrm{Int}(D_2)$ and $\psi_2\geq\psi_2(0)\equiv0$.
\end{itemize}
\end{condition}

\subsection{Well-posedness}
\label{sec:general_model_wellposedness}

\subsubsection{Well-posedness of the $X$-system}
\label{sec:general_model_wellposedness_X}

The following theorem gives the well-posedness of the $X$ process in the above system.

\bt\label{solutionX}
Under Condition \ref{X}, there exists a unique strong solution to the $X$ process in the SVI  system \eqref{eqn:OriginalsystemofFBSDEs} in the following sense:
\begin{itemize}
\item For every $t\geq0$, $X_t\in\bar{D}_1$.
 
\item For any $\varrho\in\cC(\mathbb{R}^+;\mathbb{R}^{d_1})$ and $t\geq s\geq 0$,
\begin{equation}
\label{solutionkey}
\int_s^t\langle\varrho_u-X_u,d\phi^{(1)}_u\rangle+\int_s^t\psi_1(X_u)du\leq \int_s^t\psi_1(\varrho_u)du,~~~~~a.e.,
\end{equation}
where $\phi^{(1)}$ is a continuous process of locally bounded variation, $\phi^{(1)}_0=0$.

\item For $t\in \mR_+$,
\begin{equation}
\label{eqn:Xdynamic}
X_t=x_0+\int_0^tb(s,X(s))ds+\int_{0}^t\sigma_1(s,X(s))dW_s+\int_{0}^t\sigma_2(s, X(s)) dB_s-\phi^{(1)}_t.
\end{equation}
\end{itemize}
\et

\begin{remark}\label{remark2.1}
\begin{enumerate}[(i)]
\item Note that when $\varrho=0$ in (\ref{solutionkey}), one has
$$
\int_s^t\langle X_s,d\phi^{(1)}_s\rangle\geq\int_s^t\psi_1(X_u)du.
$$
\item $\psi_1$ is locally bounded in $D_1$. Set $$
M:=\sup_{|x|\leq a}|\psi_1(x)|, \quad \varrho_u=a\frac{d\phi^{(1)}_u}{d|\phi^{(1)}|_u^r},
$$
where $|\phi^{(1)}|_u^r$ stands for the total variation of $\phi^{(1)}$ defined on an interval $[r,u]$.
Then according to equation \eqref{solutionkey},
$$
a|\phi^{(1)}|_t^s\leq\int_s^t\langle X_u,d\phi^{(1)}_u\rangle+M(t-s).
$$

\item If $(\tilde{X},\tilde{\phi}^{(1)})$ is also a solution, for any $t\geq s\geq 0$,
\begin{equation*}
\int_s^t\langle X_u-\tilde{X}_u,d\phi^{(1)}_u-d\tilde{\phi}^{(1)}_u\rangle\geq0.
\end{equation*}
\end{enumerate}
\end{remark}

We have the following lemma taken from \cite{{cepa1998problame}}. 
\bl\label{helly}
Suppose $\{k_n;n\geq1\}$ is a sequence of continuous functions $k_n: [0,T]\to\mR^{d_1}$ satisfying
$\sup_n|k_n|_0^T<\infty$ and $\|k_n-k\|_T\to0$.
Then $k$ has finite variation on $[0,T]$ and for a sequence of continuous functions $\{f_n; n\geq1\}$ satisfying $\|f_n-f\|_T\to0$ as $n\to\infty$, the following holds:
\ce
\int_s^t \langle f_n(r),dk_n(r)\rangle\to \int_s^t \langle f(r),dk(r)\rangle, \quad\mbox{as}~~n\to\infty, ~~\forall s, t\in[0,T].
\de
\el

\begin{proof} [Proof of Theorem \ref{solutionX}]
Suppose for every $T>0$ and every $n$ we are given a division of $[0,T]$:
$$
0=T_0^n<T_1^n<\cdots<T_{k_n}^n=T,
$$
with the mesh 
$$\Delta_n:=\max_{1\leq k\leq k_n}|T^n_k-T^n_{k-1}|\to0 \quad\text{as } n\to\infty.$$ 
For $t\in
(T_{k-1}^n, T_k^n]$, denote $T^n_t:=T_{k-1}^n$. Consider the following equation:
\be\label{eq02}\left\{\begin{array}{lll}
dX^n_t\in b(t,X^n(T^n_t))dt+\sigma_1(t,X^n(T^n_t))dW_t+\sigma_2(t,X^n(T^n_t)) dB_t-\partial\psi_1(X^n_t)dt, \\
X^n(0)=X^n(T_0^n)=x_0\in\bar{D}_1.
 \end{array}\right.\ee
Note that for $t\in[0,T_1^n]$, according to \cite{cepa1998problame}, there exists a unique solution to (\ref{eq02}), and we denote it by $(X^n,\phi^{(1),n})$. Applying It\^o's formula and Remark \ref{remark2.1},

\begin{align*}
|X^n_t|^2=&|x_0|^2+2\int_{0}^t\langle X^n_s, b(s,X^n(T^n_s))\rangle ds+2\int_{0}^t\langle X^n_s,\sigma_1(s,X^n(T^n_s))dW_s\rangle\\
&+2\int_{0}^t\langle X^n_s,\sigma_2(s,X^n(T^n_s))dB_s\rangle-2\int_{0}^t\langle X^n_s,d\phi^{(1),n}_s\rangle\\
&+\sum_{i=1}^2\int_{0}^t\|\sigma_i(s,X^n(T^n_s))\|^2ds\\
\leq &|x_0|^2+\int_{0}^t|X^n_s|^2ds+\int_{0}^t|b(s,X^n(T^n_s))|^2ds\\
&+\sum_{i=1}^2\int_{0}^t\|\sigma_i(s,X^n(T^n_s))\|^2ds+2M t-2a|\phi^{(1),n}|_t^{0}\\
&+2\int_{0}^t\langle X^n_s,\sigma_1(t,X^n(T^n_s))dW_s\rangle+2\int_{0}^t\langle X^n_s,\sigma_2(t,X^n(T^n_s))dB_s\rangle\\
\leq&|x_0|^2+C\int_{0}^t(1+|X^n_s|^2)ds\\
&+\int_{0}^t\bigg[b^2(s,0)+\sum_{j=1}^2\sigma_j^2(s,0)+\sum_{i=0}^2l_i^2(s)\bigg](1+\|X^n\|^2_{T^n_s})ds\\
&+2\int_{0}^t\langle X^n_s,\sigma_1(s,X^n(T^n_s))dW_s\rangle+2\int_{0}^t\langle X^n_s,\sigma_2(s,X^n(T^n_s))dB_s\rangle,
\end{align*}
from which and by using the Burkholder-Davis-Gundy (BDG) inequality and the Gr{\"o}nwall's lemma,  we have
\ce
\mathbb{E}\|X^n\|_t^2\leq C(1+\mathbb{E}|x_0|^2)\left(\int_0^T \bigg[b^2(s,0)+\sum_{j=1}^2\sigma_j^2(s,0)+\sum_{i=0}^2l_i^2(s) \bigg]ds\right),
\de
and 
\ce\begin{split}
\sup_n\mathbb{E}\sup_{t\leq T_1^n}|X^n_t|^4\leq& C(1+\mathbb{E}|x_0|^2)^2\left(\int_0^T \bigg[b^2(s,0)+\sum_{j=1}^2\sigma_j^2(s,0)+\sum_{i=0}^2l_i^2(s)\bigg]ds\right)^2.\\
\end{split}\de
Assuming \[\sup_n\mathbb{E}\sup_{t\leq T_k^n}|X^n_t|^4<\infty,\]
 then with the same arguments as above, we have
\ce
\begin{split}
&\sup_n\mathbb{E}\sup_{t\leq T_{k+1}^n}|X^n_t|^4\\
\leq& C(1+\sup_n\mathbb{E}\|X^n\|_{T_k^n}^4)\left(\int_0^T \bigg[b^2(s,0)+\sum_{j=1}^2\sigma_j^2(s,0)+\sum_{i=0}^2l_i^2(s)\bigg]ds\right)^2\\
<&\infty.
\end{split}\de
Summing up,
\be\label{momentXn}
\sup_n\mathbb{E}\sup_{t\leq T}|X^n_t|^4<\infty.
\ee

Applying It\^o's formula again, for $t\in
(T_{k-1}^n, T_k^n]$, we have
\begin{align}
\label{eqn:square_difference}
&|X^n_t-X^n_{T^n_t}|^2\nonumber\\
=&2\int_{T^n_t}^t\langle X^n_s-X^n_{T^n_t}, b(s,X^n(T^n_t))\rangle ds+2\int_{T^n_t}^t\langle X^n_s-X^n_{T^n_t}, \sigma_1(s,X^n(T^n_s))dW_s\rangle\nonumber\\
&+2\int_{T^n_t}^t\langle X^n_s-X^n_{T^n_t}, \sigma_2(s,X^n(T^n_t))dB_s\rangle-2\int_{T^n_t}^t\langle X^n_s-X^n_{T^n_t}, d\phi^{(1),n}_s\rangle\nonumber\\
&+\sum_{i=1}^2\int_{T^n_t}^t\|\sigma_i(s,X^n(T^n_s))\|^2ds\nonumber\\
\leq &\int_{T^n_t}^t|X^n_s-X^n_{T^n_t}|^2ds+\int_{T^n_t}^t|b(s,X^n(T^n_t))|^2ds+\sum_{i=1}^2\int_{T^n_t}^t\|\sigma_i(s,X^n(T^n_t))\|^2ds\nonumber\\
&+2\int_{T^n_t}^t\langle X^n_s-X^n_{T^n_t},\sigma_1(s,X^n(T^n_s))dW_s\rangle+2\int_{T^n_t}^t\langle X^n_s-X^n_{T^n_t},\sigma_2(s,X^n(T^n_s))dB_s\rangle\\
&-2\int_{T^n_t}^t\langle X^n_s-X^n_{T^n_t}, d\phi^{(1),n}_s\rangle.\nonumber
\end{align}

For $\epsilon>0$ and $R>0$, set
\be\label{Ae}\begin{split}
A_{\epsilon,R}:=&\{x\in\mR^{d_1}:\forall x'\notin\bar{D}_1,|x-x'|\geq\epsilon~~\mbox{and}~~|x-a_0|\leq R\},
\end{split}
\ee
where $a_0\in\mathrm{Int}(D_1)$
 such that $A_{\epsilon,R}\neq\emptyset$~ for every $R>0$ and $\epsilon<\epsilon_0$ for some $\epsilon_0>0$.
Then $A_{\epsilon,R}$ is a convex compact subset of $\mathrm{Int}(D_1)$.
Set $$f_R(\epsilon):=\sup\{|x'|:x'\in \partial\psi_1(x), x\in A_{\epsilon,R}\},$$
and according to the local boundedness of $\partial\psi_1$ on $\mathrm{Int}(D_1)$, $|f_R(\epsilon)|<+\infty$.
 Let \ce
g_R(\delta):=\inf\{\epsilon\in(0,\epsilon_0):f_R(\epsilon)\leq\delta^{-1/2}\},\quad
\delta>0.\de
Let
$\delta_R>0$ such that $\delta_R+g_R(\delta_R)<\epsilon_0$. Fix
$R>0$ and $\delta\in(0,\delta_R\wedge 1]$. Since
$$\delta_R+g_R(\delta_R)<\epsilon_0, \quad A_{\delta+g_R(\delta),R}\neq\emptyset,$$
we have 
\begin{equation} \label{eqn:f_R_bound} f_R(\delta+g_R(\delta))\leq
\delta^{-1/2}.\end{equation} 

For $0\leq t-s\leq\delta$, denote 
$\xi^{n,\delta,R}$ as the projection of $X^n_s$ on
$A_{\delta+g_R(\delta),R}$. 
Then on the set $\{\|X^n\|_T\leq R\}$, we have
\ce
|X^n_{T^n_t}-\xi^{n,\delta,R}|\leq \delta+g_R(\delta),
\de
which yields
$$\int_{s}^t\langle X^n_{s}-\xi^{n,\delta,R}, d\phi^{(1),n}_r\rangle\leq (\delta+g_R(\delta))|\phi^{(1),n}|_T^0,$$
and
\begin{align*}
\int_{s}^t\langle\xi^{n,\delta,R}-X^n_r,d\phi^{(1),n}_r\rangle
\leq &\int_{s}^t\langle\xi^{n,\delta,R}-X^n_r,\eta^{n,\delta,R}\rangle dr\\
\leq &2R(t-s)f_R(\delta+g_R(\delta))\\
\leq &2\delta^{1/2}R,
\end{align*}
where the first inequality follows by equation \eqref{solutionkey} with $\eta^{n,\delta,R}\in \partial\psi_1(\xi^{n,\delta,R})$, the second inequality follows by the boundedness of $\xi^{n,\delta,R}$ and the definitions of $\xi^{n,\delta,R}$ and $f_R(\delta+g_R(\delta))$, and the third inequality follows by equation \eqref{eqn:f_R_bound}.
Therefore, on the set $\{\|X^n\|_T\leq R\}$,
\begin{equation}\begin{split}
\label{eqn:inner_product_X_phi}
-\int_{s}^t\langle X^n_r-X^n_{s}, d\phi^{(1),n}_r\rangle
=&\int_{s}^t\langle X^n_{s}-\xi^{n,\delta,R}, d\phi^{(1),n}_r\rangle+\int_{s}^t\langle\xi^{n,\delta,R}-X^n_r,d\phi^{(1),n}_r\rangle\\
\leq & (\delta+g_R(\delta))|\phi^{(1),n}|_T^0+2\delta^{1/2}R.
\end{split}\end{equation}

Define the stopping time 
$$\tau_n(R):=\inf\{s;|X^n_s|>R\}.$$ For $t\leq \tau_n(R)\wedge T$, 
 plugging the result of \eqref{eqn:inner_product_X_phi} in \eqref{eqn:square_difference}, we have
\ce\begin{split}
&|X^n_t-X^n_{T^n_t}|^2\\
\leq &\int_{T^n_t}^t|X^n_s-X^n_{T^n_t}|^2ds+\int_{T^n_t}^t|b(s,X^n(T^n_t))|^2ds+\sum_{i=1}^2\int_{T^n_t}^t\|\sigma_i(s,X^n(T^n_t))\|^2ds\\
&+2\int_{T^n_t}^t\langle X^n_s-X^n_{T^n_t},\sigma_1(s,X^n(T^n_s))dW_s\rangle+2\int_{T^n_t}^t\langle X^n_s-X^n_{T^n_t},\sigma_2(s,X^n(T^n_s))dB_s\rangle\\
&+2(\Delta_n+g_R(\Delta_n))|\phi^{(1),n}|_T^0+4R\Delta_n^{1/2}.
\end{split}\de
Taking supremum and then expectation, we have
\be\label{Deltan}\begin{split}
\mathbb{E}\sup_{t\leq T\wedge \tau_n(R)}|X^n_t-X^n_{T^n_t}|^2\leq& C\Delta_n^{1/2}(1+\mathbb{E}\|X^n\|_T^2)+\mathbb{E}|\phi^{(1),n}|_T^0(\Delta_n+g_R(\Delta_n))\\
&+\max_k\int_{T_k^n}^{T_{k+1}^n}\bigg[b^2(s,0)+\sum_{j=1}^2\sigma_j^2(s,0)+\sum_{i=0}^2l_i^2(s)\bigg]ds,
\end{split}\ee
which together with (\ref{momentXn}) implies that
\ce\begin{split}
\mathbb{E}\sup_{t\leq T}|X^n_t-X^n_{T^n_t}|^2\leq&\mathbb{E}\sup_{t\leq T}|X^n_t-X^n_{T^n_t}|^2(\mathbbm{1}_{\{T<\tau_n(R)\}}+\mathbbm{1}_{\{T\geq\tau_n(R)\}})\\
\leq&\mathbb{E}\sup_{t\leq T\wedge \tau_n(R)}|X^n_t-X^n_{T^n_t}|^2+\mathbb{E}\|X^n\|_T^2 \mathbbm{1}_{\{T\geq\tau_n(R)\}}\\
\to &0, \quad \mbox{by letting} ~~~n\to\infty ~~~~\mbox{and then} ~~~~R\to\infty.
\end{split}\de

Furthermore, by Condition \ref{X} which implies that 
$$\int_{0}^t \big\langle X^n_s-X^m_s,b(s,X^n(s))-b(s,X^m(s))\big\rangle ds\leq 0,$$
and by Remark \ref{remark2.1} which implies that
$$\int_{0}^t\big\langle X^n_s-X^m_s,d(\phi^{(1),n}_s-\phi^{(1),m}_s)\big\rangle\geq 0,$$
we have
\begin{align*}
&|X^n_t-X^m_t|^2\\
=&2\int_{0}^t\big\langle X^n_s-X^m_s, b(s,X^n(T^n_s))-b(s,X^m(T^m_s))\big\rangle ds\\
&+2\int_{0}^t \big\langle X^n_s-X^m_s, \sigma_1(s,X^n(T^n_s))-\sigma_1(s,X^m(T^m_s))\big\rangle dW_s\\
&+2\int_{0}^t \big\langle X^n_s-X^m_s,\sigma_2(s,X^n(T^n_s))-\sigma_2(s,X^m(T^m_s))dB_s\big\rangle\\
&-2\int_{0}^t \big\langle X^n_s-X^m_s ,d(\phi^{(1),n}_s-\phi^{(1),m}_s)\big\rangle\\
&+\sum_{i=1}^2\int_{0}^t \|\sigma_i(s,X^n(T^n_s))-\sigma_i(s,X^m(T^m_s))\|^2ds\\
\leq&2\int_{0}^t \big\langle X^n_s-X^m_s,b(s,X^n(T^n_s))-b(s,X^n(s))\big\rangle ds\\
&+2\int_{0}^t \big\langle X^n_s-X^m_s,b(s,X^m(T^m_s))-b(s,X^m(s))\big\rangle ds\\
&+2\int_{0}^t \big\langle X^n_s-X^m_s,\sigma_1(s,X^n(T^n_s))-\sigma_1(s,X^m(T^m_s))\big\rangle dW_s\\
&+2\int_{0}^t \big\langle X^n_s-X^m_s,\sigma_2(s,X^n(T^n_s))-\sigma_2(s,X^m(T^m_s))\big\rangle dB_s\\
&+\sum_{i=1}^2\int_{0}^tl_i^2(s)\|X^m(T^m_\cdot)-X^n(T^n_\cdot)\|_s^2ds\\
\leq&2\int_{0}^tl_0(s)|X^n_s-X^m_s|\big(|X^n(T^n_s)-X^n(s)|^{\frac12+\alpha}+|X^m(T^m_s)-X^m(s)|^{\frac12+\alpha}\big) ds\\
&+2\int_{0}^t \big\langle X^n_s-X^m_s, \sigma_1(s,X^n(T^n_s))-\sigma_1(s,X^m(T^m_s))\big\rangle dW_s\\
&+2\int_{0}^t \big\langle X^n_s-X^m_s, \sigma_2(s,X^n(T^n_s))-\sigma_2(s,X^m(T^m_s))\big\rangle dB_s\\
&+C\sum_{i=1}^2\int_{0}^tl_i^2(s)\big(\|X^m(T^m_\cdot)-X^m(\cdot)\|_s^2+\|X^n(T^n_\cdot)-X^n(\cdot)\|_s^2\big)ds\\
&+C\sum_{i=1}^2\int_{0}^tl_i^2(s)\|X^m-X^n\|_s^2 ds.
\end{align*}

Define the stopping time $$\tau_m(R):=\inf\{s;|X^m_s|>R\}.$$ On one hand, by the BDG inequality and equation \eqref{Deltan}, we get
\ce\begin{split}
&\mathbb{E}\sup_{t\leq T}|X^n_t-X^m_t|^2\mathbbm{1}_{\{T<\tau_m(R)\wedge\tau_n(R)\}}\\
=&\mathbb{E}\sup_{t\leq T\wedge\tau_m(R)\wedge\tau_n(R)}|X^n_t-X^m_t|^2\\
\leq&C_T\big(h_R(\Delta_m)+h_R(\Delta_n)\big)+C\sum_{i=1}^2\int_{0}^tl_i^2(s)\mathbb{E}\|X^m-X^n\|_s^2ds,
\end{split}\de
where $h_R(\Delta_k)\to0$ as $k\to\infty$. 
On the other hand, by H\"older's inequality and equation \eqref{momentXn},
\ce\begin{split}
\mathbb{E}\sup_{t\leq T}|X^n_t-X^m_t|^2\mathbbm{1}_{\{T\geq\tau_m(R)\wedge\tau_n(R)\}}
\leq & \left[ \mathbb{E}\sup_{t\leq T}|X^n_t-X^m_t|^4 \cdot \mathbb{E}\mathbbm{1}_{\{T\geq\tau_m(R)\wedge\tau_n(R)\}} \right]^{\frac{1}{2}}\\
\leq & \left[2\sup_n \mathbb{E}\sup_{t\leq T}|X^n_t|^4 \cdot \mathbb{P}(T\geq\tau_m(R)\wedge\tau_n(R)) \right]^{\frac{1}{2}}\\
\leq & \left[C \cdot \mathbb{P} \left(\sup_{t\leq T}|X^n_t| \vee \sup_{t\leq T}|X^m_t| >R\right) \right]^{\frac{1}{2}}\\
\leq & \left[\frac{C}{R^2} \sup_n \mathbb{E}\sup_{t\leq T}|X^n_t|^2 \right]^{\frac{1}{2}}\\
\leq & \frac{C}{R^2}.
\end{split}\de
Hence,
\ce\begin{split}
\mathbb{E}\sup_{t\leq T}|X^n_t-X^m_t|^2=&\mathbb{E}\bigg[\sup_{t\leq T}|X^n_t-X^m_t|^2(\mathbbm{1}_{\{T<\tau_m(R)\wedge\tau_n(R)\}}+\mathbbm{1}_{\{T\geq\tau_m(R)\wedge\tau_n(R)\}})\bigg]\\
\leq& C_T\big(h_R(\Delta_m)+h_R(\Delta_n)\big)+\frac{C}{R^2}\\
\to& 0, ~~~~\mbox{as}~~~m, n\to\infty~~~~\mbox{and then}~~~~R\to\infty,
\end{split}\de
and moreover by equation \eqref{eqn:Xdynamic},
$$
\lim_{m,n\to\infty}\mathbb{E}\|\phi^{(1),m}-\phi^{(1),n}\|_T\to0.
$$

Hence, $\{X_n,\phi^{(1),n}\}_n$ is a Cauchy sequence and by the completeness of the space of processes with respect to the uniform convergence, there exists a pair of continuous processes $(X,\phi^{(1)})$ satisfying that for any $\epsilon>0$,
\ce
\mathbb{E}\sup_{t\leq T}|X^n_t-X_t|^2\rightarrow0,~~~~\mathbb{E}\sup_{t\leq T}\big|\phi^{(1),n}_t-\phi_t^{(1)}\big|^2\rightarrow0.
\de
Then by Lemma \ref{helly}, we have that $\phi^{(1)}$ is of locally finite variations and equation (\ref{solutionkey}) holds.
Furthermore, by the continuity of $b$ and $\sigma$, we have
\begin{align*}
&\mathbb{E}\sup_{t\leq T}\left|\int_0^{t}\sigma_1(s,X^n(T^n_s))d{W}_s-\int_0^{t}\sigma_1(s,{X}(s))d{W}_s\right|^2\to0,\\
&\mathbb{E}\sup_{t\leq T}\left|\int_0^{t}\sigma_2(s,X^n(T^n_s))d{B}_s-\int_0^{t}\sigma_2(s,{X}(s))d{B}_s\right|^2\to0,\\
&\mathbb{E}\sup_{t\leq T}\left|\int_0^{t}b(s,{X}^n(T^n_s))ds-\int_0^{t}b(s,{X}(s))ds\right|^2\to0.
\end{align*}

Suppose $(\bar{X},\bar{\phi}^{(1)})$ is also a solution. It\^o's formula with Remark \ref{remark2.1} yields
\begin{align*}
&|X_t-\bar{X}_t|^2\\
=&2\int_0^t\big\langle X_s-\bar{X}_s, b(s,X(s))-b(s,\bar{X}(s))\big\rangle ds\\
&+2\int_0^t\big\langle X_s-\bar{X}_s,\sigma_1(s,X(s))-\sigma_1(s,\bar{X}(s))\big\rangle dW_s\\
&+2\int_0^t\big\langle X_s-\bar{X}_s, \sigma_2(s,X(s))-\sigma_2(s,\bar{X}(s))\big\rangle dB_s\\
&-2\int_0^t\big\langle X_s-\bar{X}_s,\big(d\phi^{(1)}_s-d\bar{\phi}^{1}_s\big)\big\rangle+\sum_{j=1}^2\int_0^t\|\sigma_j(s,X(s))-\sigma_j(s,\bar{X}(s))\|^2ds\\
\leq&2\int_0^t\big\langle X_s-\bar{X}_s, \sigma_1(s,X(s))-\sigma_1(s,\bar{X}(s)\big\rangle dW_s\\
&+2\int_0^t\big\langle X_s-\bar{X}_s,\sigma_2(s,X(s))-\sigma_2(s,\bar{X}(s))\big\rangle dB_s+\sum_{j=1}^2\int_0^tl_j^2(s)\|X-\bar{X}\|_s^2ds,
\end{align*}
from which we could get
\ce
\mathbb{E}\|X-\bar{X}\|_T^2\leq C\sum_{j=1}^2\int_0^Tl_j^2(s)\mathbb{E}\|X-\bar{X}\|_s^2ds,
\de
and the uniqueness follows by Gr{\"o}nwall's inequality.
\end{proof}

\subsubsection{Well-posedness of the $Y$-system}
\label{sec:general_model_wellposedness_Y}

\begin{remark}\label{remark2.1Y}
Analogous to Theorem \ref{solutionX} and Remark \ref{remark2.1Y}, one can show that 
\begin{itemize}
\item For any $\varrho\in\mathcal{C}(\mathbb{R}^+;\mathbb{R}^{d_2})$ and $t\geq s\geq 0$,
\begin{equation}
\label{solutionkeyY}
\int_s^t \big\langle\varrho_u-Y_u,d\phi^{(2)}_u\big\rangle+\int_s^t\psi_2(Y_u)du\leq \int_s^t\psi_2(\varrho_u)du,~~~~~a.e.,
\end{equation}
where $\phi^{(2)}$ is a continuous process of locally bounded variation satisfying that $\phi^{(2)}_0=0$.

\item  If $(Y,\phi^{(2)})$ and $(\tilde{Y},\tilde{\phi}^{(2)})$ are two solutions, then for any $t\geq s\geq 0$,
\begin{equation*}
\int_s^t\big\langle Y_u-\tilde{Y}_u,d\phi^{(2)}_u-d\tilde{\phi}^{(2)}_u\big\rangle\geq 0.
\end{equation*}
\end{itemize}
\end{remark}

\begin{proposition}\label{wellY}
Under Conditions \ref{X} and \ref{Y}, there exists a unique strong solution to the $Y$ process in the SVI  system \eqref{eqn:OriginalsystemofFBSDEs}.
\end{proposition}

\begin{proof}
Suppose $Z$ is an adapted process satisfying
\ce
\mathbb{E}\|Z\|_T^4<\infty.
\de
Then according to the deterministic result (see \cite{cepa1998problame}), there exists a unique solution $(Y,\phi^{(2)})$ to the following SVI:
\be\label{YZ}
Y_t\in y_0+\int_0^t \alpha(s,X(s),Z(s),q_s)ds+\int_0^t\beta(s,X(s),Z(s),q_s) dB_s-\int_0^t\partial\psi_2(Y_s)ds.
\ee

Note that similar to (\ref{momentXn}), we have 
\be\label{Xmoment}
\mathbb{E}\|X\|_T^4<\infty.
\ee

Denote $$\tau_R^1:=\inf\{s;|X_s|\vee|Z_s|>R\}.$$
Then for all $R>0$, $\tau_R^1$ is a stopping time and $\tau_R^1\uparrow\infty$ as $R\uparrow\infty$. By It\^o's formula and with arguments similar to the previous section, for any $t<\tau_R^1$,
\begin{align*}
|Y_t|^2\leq&|y_0|^2+2\int_0^t \big\langle Y_s, \alpha(s,X(s),Z(s),q_s)\big\rangle ds-2\int_0^t \big\langle Y_s, d\phi^{(2)}_s \big\rangle\\
&+\int_0^t\|\beta(s,X(s),Z(s),q_s)\|^2 ds+2\int_0^t \big\langle Y_s,\beta(s,X(s),Z(s),q_s)dB_s\big\rangle\\
\leq&|y_0|^2+\int_0^t|Y_s|^2ds+\int_0^t \left(L_R^2(s)\|Z\|_s^2+|\alpha(s,X(s),0,q_s)|^2\right)ds\\
&+\lambda_2^2\int_0^t(L_R^2(s)\|Z\|_s^2+\|\beta(s,X(s),0,q_s)\|^2)ds
+2M t-2a|\phi^{(2)}|_t^{0}\\&+2\int_0^t \big\langle Y_s,  \beta(s,X(s),Z(s),q_s) dB_s\big\rangle,
\end{align*}
where in the last inequality we used equation \eqref{solutionkeyY} and the mean value theorem.
Hence,
\begin{align*}
&\mathbb{E}\sup_{t\leq T\wedge\tau_R^1}|Y_t|^4\\
\leq&C\mathbb E(1+|y_0|^4)+C\mathbb E\int_0^{T\wedge\tau_R^1}|Y_s|^4ds\\
&+C\mathbb{E}\left[\int_0^{T\wedge\tau_R^1}\left(L_R^2(s)\|Z\|_s^2+|\alpha(s,X(s),0,q_s)|^2\right)ds\right]^2\\
&+C\mathbb{E}\lambda_2^4\left[\int_0^{T\wedge\tau_R^1}\left(L_R^2(s)\|Z\|_s^2+\|\beta(s,X(s),0,q_s)\|^2\right)ds\right]^2\\
&+C\mathbb{E}\int_0^{T\wedge\tau_R^1}|Y_s|^2\cdot \|\beta(s,X(s),Z(s),q_s)\|^2 ds\\
\leq&C\mathbb E(1+|y_0|^4)+C\mathbb E\int_0^{T\wedge\tau_R^1}|Y_s|^4ds
+\frac12\mathbb{E}\sup_{t\leq T\wedge\tau_R^1}|Y_t|^4\\
&+C\mathbb{E}\|Z\|_T^4\bigg[\int_0^{T\wedge\tau_R^1}\bigg(L_R^2(s)\\
&\quad\quad\quad\quad\quad\quad\quad\quad\quad+\sup_{\|x\|_T\leq R,\lambda_1\leq\|y\|\leq\lambda_2}
(|\alpha(s,x,0,y)|^2+\|\beta(s,x,0,y)\|^2)\bigg)ds\bigg]^2,
\end{align*}
and thus by the Gr{\"o}nwall's lemma
$$\mathbb{E}\sup_{t\leq T\wedge\tau_R^1}|Y_t|^4\leq C (1+\mathbb{E}\|Z\|_T^4).$$
Furthermore,
\begin{equation}
\begin{split}
\label{eqn:supY_tail}
\mathbb{P}(\|Y\|_T>M)
=&\mathbb P(\|Y\|_T>M,T<\tau_R^1)+\mathbb{P}(\|Y\|_T>M,T\geq\tau_R^1)\\
\leq &\mathbb P(\|Y\|_{T\wedge\tau_R^1}>M)+\mathbb{P}(T\geq\tau_R^1)\\
\leq&\frac{\mathbb{E}\|Y\|_{T\wedge\tau_R^1}^4}{M^4}+\mathbb{P}(T\geq\tau_R^1)\\
\to&0, \quad\mbox{by letting} ~~~M\to\infty ~~\mbox{and then} ~~~R\to\infty.
\end{split}
\end{equation}

Now we are going to show that the map $Z\to (Y_{\cdot\wedge\tau_R^1},\phi_{\cdot\wedge\tau_R^1}^{(2)})$ is a contraction. Suppose $\bar{Z}$ is also an adapted process such that
$$
\mathbb{E}\|\bar{Z}\|^4_T<\infty,
$$
and $(\bar{Y},\bar{\phi}^{(2)})$ is the unique solution to equation \eqref{YZ} with $\bar{Z}$ in place of $Z$. Define
$$\tau_R^1:=\inf\{s;|X_s|\vee|Z_s|\vee|\bar{Z}_s|>R\}.$$ 
Remark \ref{remark2.1Y} implies that
$$\int_{0}^t \big\langle Y_s-\bar{Y}_s, d(\phi^{(2)}_s-\bar{\phi}^{(2)}_s)\big\rangle \geq 0,$$
and then by It\^o formula and Condition \ref{Y}
\ce
\begin{split}
&|Y_{t\wedge\tau_R^1}-\bar{Y}_{t\wedge\tau_R^1}|^2\\
\leq&\int_0^{t\wedge\tau_R^1}|Y_s-\bar{Y}_s|^2ds+\int_0^{t\wedge\tau_R^1}|\alpha(s,X(s),Z(s),q_s)-\alpha(s,X(s),\bar{Z}(s),q_s)|^2ds\\
&+\int_0^{t\wedge\tau_R^1}\|\beta(s,X(s),Z(s),q_s)-\beta(s,X(s),\bar{Z}(s),q_s)\|^2ds\\
&+2\int_0^{t\wedge\tau_R^1} \big\langle(Y_s-\bar{Y}_s), \big(\beta(s,X(s),Z(s),q_s)-\beta(s,X(s),\bar{Z}(s),q_s)\big) dB_s\big\rangle\\
\leq&\int_0^{t\wedge\tau_R^1}|Y_s-\bar{Y}_s|^2ds+(1+\lambda_2^2)\int_0^{t\wedge\tau_R^1}L_R^2(s)\|Z-\bar{Z}\|_s^2ds\\
&+2\int_0^{t\wedge\tau_R^1}\big\langle(Y_s-\bar{Y}_s), \big(\beta(s,X(s),Z(s),q_s)-\beta(s,X(s),\bar{Z}(s),q_s)\big) dB_s\big\rangle.
\end{split}\de
Set $l_t:=\int_0^t L_R^2(s)ds$. Taking supremum and expectation of the above equation yields 
\ce
\begin{split}
\mathbb{E}\sup_{t\leq T}|Y_{t\wedge\tau_R^1}-\bar{Y}_{t\wedge\tau_R^1}|^2
\leq&C(\lambda_2,T)\mathbb{E}\int_0^{T\wedge\tau_R^1}L_R^2(s)\|Z-\bar{Z}\|_s^2ds\\
\leq&C(\lambda_2,T)\left(\int_0^{T}L_R^2(s)e^{rl_s}ds\right) \cdot \left(\sup_{t\leq T}e^{-rl_t}\mathbb E\|Z-\bar{Z}\|_t^2\right)\\
=&\frac{C(\lambda_2,T)}{r}e^{rl_T}\sup_{t\leq T}e^{-rl_t}\mathbb{E}\|Z-\bar{Z}\|_t^2.
\end{split}\de
Taking $r=2C(\lambda_2,T)$ gives
\ce
\sup_{t\leq T}e^{-rl_t}\mathbb E\|Y-\bar{Y}\|_t^2\leq \frac12\sup_{t\leq T}e^{-rl_t}\mathbb{E}\|Z-\bar{Z}\|_t^2.
\de
Let $Y^{(0)}\equiv y$ and for $n\geq1$, denote $(Y^{n},\phi^{(2),n})$ as the solution to equation (\ref{YZ}) with $Z$ replaced by $Y^{n-1}$. Then for any $\delta>0$,
\ce\begin{split}
\mathbb{P}(\|Y^{n}-Y^{n-1}\|_T>\delta)
\leq&\mathbb{P}(\|Y^{n}-Y^{n-1}\|_T>\delta,T<\tau_R^1)+\mathbb{P}(T\geq \tau_R^1)\\
\leq&\frac{e^{2rl_T}}{\delta^2}e^{-rl_T}\mathbb{e}\|Y^{n}-Y^{n-1}\|_{T\wedge\tau_R^1}^2+\mathbb{P}(T\geq \tau_R^1)\\
\leq&\frac{e^{2rl_T}}{\delta^2}\left(\frac12 \right)^{n-1}\mathbb{E}\|Y^{1}\|_{T\wedge\tau_R^1}^2+\mathbb{P}(T\geq \tau_R^1)\\
\to&0, \quad\quad\mbox{by letting}~n\to\infty ~\mbox{and then}~R\to\infty,\\
\end{split}
\de
which, by the $Y$ dynamic, yields
$$\mathbb P(\|\phi^{(2),n}-\phi^{(2),n-1}\|_T>\delta)\to 0, \quad\quad\mbox{by letting}~n\to\infty.$$
Thus, by completeness there exists a unique pair of processes $(Y,\phi^{(2)})$ such that
\begin{align*}
&\mathbb P(\|Y^{n}-Y\|_T>\delta)\to0,\quad \mathbb P(\|\phi^{(2),n}-\phi^{(2)}\|_T>\delta)\to0, &\quad\mbox{by letting}~n\to\infty.
\end{align*}
By equation \eqref{eqn:supY_tail} we have that 
$$\mathbb{P}(\|Y^{n}\|_T>M)\to0, \quad  \mathbb{P}(|\phi^{(2)}|_T>M)\to0,\quad\mbox{as}~M\to\infty,$$ from which we get 
$$\mathbb P(\|Y\|_T>M)\to0, \quad\mathbb P(|\phi^{(2)}|_T^0>M)\to0, \quad \mbox{as}~M\to\infty.$$
Applying Lemma \ref{helly}, for any $a\in\bar{D}_2$ and $t\geq s\geq r$,
$$
\int_s^t(a-Y_r)d\phi^{(2)}_r+\int_s^t\psi_2(Y_r)dr\leq (t-s)\psi_2(a),~~~~~a.e..
$$
Hence we have proved that $(Y,\phi^{(2)})$ is a solution of the $Y$ process in the SVI  system \eqref{eqn:OriginalsystemofFBSDEs}.

To prove the uniqueness, we first suppose $(\tilde{Y},\tilde{\phi}^{(2)})$ is also a solution. Denote
$$\tau_R:=\inf\{s;|X_s|\vee|Y_s|\vee|\tilde{Y}_s|>R\}.$$ Applying It\^o's formula, for $t<\tau_R$, yields
\ce
\begin{split}
|Y_t-\tilde{Y}_t|^2\leq&2\int_0^t \big\langle Y_s-\tilde{Y}_s,\big[\alpha(s,X(s),Y(s),q_s)-\alpha(s,X(s),\tilde{Y}(s),q_s)\big]\big\rangle ds\\
&+2\int_0^t \big\langle Y_s-\tilde{Y}_s,\big[\beta(s,X(s),Y(s),q_s)-\beta(s,X(s),\tilde{Y}(s),q_s)\big] dB_s\big\rangle\\
&+\int_0^t\|\beta(s,X(s),Y(s),q_s)-\beta(s,X(s),\tilde{Y}(s),q_s)\|^2ds.
\end{split}\de
Then taking expectations yields
\ce
\begin{split}
\mathbb{E}\sup_{t\leq T\wedge\tau_R}|Y_t-\tilde{Y}_t|^2\leq C\mathbb{E}\int_0^{T\wedge\tau_R}|Y_s-\tilde{Y}_s|^2ds+C\mathbb{E}\int_0^{T\wedge\tau_R}\|Y-\tilde{Y}\|_s^2ds,
\end{split}\de
from which we have
\ce
\mathbb{E}\sup_{t\leq T\wedge\tau_R}|Y_t-\tilde{Y}_t|^2=0,
\de
and furthermore
\ce
\mathbb{P}\left(\sup_{t\leq T\wedge\tau_R}|Y_t-\tilde{Y}_t|>0\right)=0.
\de
\end{proof}

\subsection{Asymptotic Analysis}
\label{sec:general_model_asymptotic}

We now study the stability of the SVI  system \eqref{eqn:OriginalsystemofFBSDEs} by investigating its perturbed version with a small positive parameter $\epsilon$
\begin{equation}
\label{eqn:SVI_P}
\left\{\begin{array}{lll}
 X_t^{\varepsilon}\in & x_0+\int_0^tb^{\varepsilon}(s,X^{\varepsilon}(s),\varepsilon)ds+\int_0^t\sigma_1^{\varepsilon}(s,X^{\varepsilon}(s),\varepsilon)dW_s+\int_0^t\sigma_2^{\varepsilon}(s, X^{\varepsilon}(s), \varepsilon) dB_s\\
 &-\int_0^t\partial\psi_1(X^{\varepsilon}_s)ds,\\
       \\
Y_t^{\varepsilon}\in & y_0+\int_0^t \alpha(s,X^{\varepsilon}(s),Y^{\varepsilon}(s), q_s)ds+\int_0^t\beta(s,X^{\varepsilon}(s),Y^{\varepsilon}(s), q_s)dB_s\\
&-\int_0^t\partial\psi_2(Y^{\varepsilon}_s)ds,
             \end{array}\right.
\end{equation}
where 
\begin{equation}
\label{eqn:barX}
\lim_{\varepsilon \rightarrow 0} b^{\varepsilon}(t, x, \varepsilon)=b(t, x), \quad \lim_{\varepsilon \rightarrow 0} \sigma_i^{\varepsilon}(t, x, \varepsilon)=\sigma_i(t, x), \quad i=1, 2.
\end{equation}

\begin{condition}\label{Cond:asymptotic_X}
Suppose that $b^{\varepsilon}(t, x, \varepsilon)$ and $\sigma_j^{\varepsilon}(t, x, \varepsilon)$ for $j=1,2$ are continuous in $t$ uniformly in $\varepsilon$, and satisfy
\begin{align*}
&\big\langle b^{\varepsilon}(t, x(t), \varepsilon)-b^{\varepsilon}(t, x'(t),\varepsilon), x(t)-x'(t)\big\rangle\leq 0,\quad &\forall x, x'\in\cC(\mathbb R^+;\mathbb R^{d_1}),\\
&|b^{\varepsilon}(t, x(t), \varepsilon)-b^{\varepsilon}(t, x'(t), \varepsilon)|\leq l_0(t) \|x-x'\|_t^{1/2+\alpha}, \quad &\text{for some } \alpha\in[0,1/2],\\
&\|\sigma_i^{\varepsilon}(t, x(t), \varepsilon)-\sigma_i^{\varepsilon}(t, x'(t),\varepsilon)\|\leq l_i(t)\|x-x'\|_t, \quad & i=1, 2,
\end{align*}
where $l_i(t)$ for $i=0,1,2$ are functions of $t$ satisfying that $l_i(\cdot)\in L^2([0,T])$. 
\end{condition}

\subsubsection{Asymptotic analysis of the $X$ system}
\label{sec:general_model_asymptotic_X}

In the following, we give the convergence result regarding the $X_t^{\varepsilon}$ process in the perturbed system \eqref{eqn:SVI_P} as $\varepsilon$ goes  to $0$.
\bt\label{Xpathasymp}
As $\varepsilon \rightarrow 0$, under Conditions \ref{X} and \ref{Cond:asymptotic_X}, we have
\begin{equation}
\mathbb{E}\sup_{t\in [0,T]}|X_t^{\varepsilon}-X_t|^2\rightarrow 0.
\end{equation}
 \et

\begin{proof}
By applying It\^o's formula,
\begin{align*}
|X_t^{\varepsilon}-X_t|^2=&2\int_0^t\big\langle X_s^{\varepsilon}-X_s,b^{\varepsilon}(s,X^{\varepsilon}(s),\varepsilon)-b(s,X(s))\big\rangle ds\\
&+\sum_{i=1}^2\int_0^t\|\sigma_i^{\varepsilon}(s,X^{\varepsilon}(s),\varepsilon)-\sigma_i(s,X(s))\|^2ds\\
&+2\int_0^t \big\langle X_s^{\varepsilon}-X_s,\big(\sigma_1^{\varepsilon}(s,X^{\varepsilon}(s),\varepsilon)-\sigma_1(s,X(s))\big)dW_s\big\rangle\\
&+2\int_0^t \big\langle X_s^{\varepsilon}-X_s,\big(\sigma_2^{\varepsilon}(s,X^{\varepsilon}(s),\varepsilon)-\sigma_2(s,X(s))\big)dB_s\big\rangle\\
&-2\int_0^t \big\langle X_s^{\varepsilon}-X_s ,d\phi^{(1),\varepsilon}_s-d\phi^{(1)}_s\big\rangle\\
\leq&C\int_0^t\big(1+l_1^2(s)+l_2^2(s)\big)\|X^{\varepsilon}-X\|_s^2ds\\
&+\int_0^t|b^{\varepsilon}(s,\bar{X}(s),\varepsilon)-b(s,X(s))|^2ds\\
&+\sum_{i=1}^2\int_0^t\|\sigma_i^{\varepsilon}(s,X(s),\varepsilon)-\sigma_i(s,X(s))\|^2ds\\
&+2\int_0^t \big\langle X_s^{\varepsilon}-X_s,\big(\sigma_1^{\varepsilon}(s,X^{\varepsilon}(s),\varepsilon)-\sigma_1(s,X(s))\big)dW_s\big\rangle\\
&+2\int_0^t \big\langle X_s^{\varepsilon}-X_s,\big(\sigma_2^{\varepsilon}(s,X^{\varepsilon}(s),\varepsilon)-\sigma_2(s,X(s))\big)dB_s\big\rangle,
\end{align*}
which implies that
\ce
\begin{split}
\mathbb{E}\|X^{\varepsilon}-X\|_T^2\leq&C\mathbb{E}\int_0^T\big(1+l_1^2(s)+l_2^2(s)\big)\|X^{\varepsilon}-X\|_s^2ds\\
&+C\mathbb{E}\int_0^T|b^{\varepsilon}(s,X(s),\varepsilon)-b(s,X(s))|^2ds\\
&+\sum_{i=1}^2\mathbb{E}\int_0^T\|\sigma_i^{\varepsilon}(s,X(s),\varepsilon)-\sigma_i(s,X(s))\|^2ds.
\end{split}\de
The Gr{\"o}nwall's lemma yields that
\ce
\begin{split}
\mathbb{E}\|X^{\varepsilon}-X\|_T^2\leq&C\mathbb{E}\int_0^T|b^{\varepsilon}(s,X(s),\varepsilon)-b(s,X(s))|^2ds\\
&+C\sum_{i=1}^2\mathbb{E}\int_0^T\|\sigma_i^{\varepsilon}(s,X(s),\varepsilon)-\sigma_i(s,X(s))\|^2ds.\end{split}\de
Now it follows from (\ref{Xmoment}) and (\ref{eqn:barX}) that 
\ce\begin{split}
\mathbb{E}\|X^{\varepsilon}-X\|_T^2\to  0,\quad\quad\quad\mbox{as}\;\; \varepsilon\to0.
\end{split}\de
\end{proof}

\subsubsection{Asymptotic analysis of the $Y$ system}
\label{sec:general_model_asymptotic_Y}

In the following, we give the convergence result regarding the $Y_t^{\varepsilon}$ process in the perturbed system \eqref{eqn:SVI_P} as $\varepsilon$ goes  to $0$.
\bt \label{thm:Ycvg}
Under Conditions \ref{X}, \ref{Y}, and \ref{Cond:asymptotic_X}, as $\varepsilon \rightarrow 0$, for any $\eta>0$, we have
\begin{equation}
\mathbb{P}\left(\sup_{t\in [0,T]}|Y_t^{\varepsilon}-Y_t|>\eta\right)\rightarrow 0.
\end{equation}
\et

\begin{proof}
We firstly define stopping time $\tau$ as
\begin{equation}
\label{eqn:tau}
\tau=\inf\{s: |X_s^{\varepsilon}|>R\}.
\end{equation}
Then with analysis analogous to Proposition \ref{wellY}, we have
\begin{equation*}
\begin{split}
&\mathbb{E}\sup_{t\in [0,T]}|Y_{t \wedge \tau}^{\varepsilon}|^2
\\
\leq & |y_0|^2+C\int_0^{T \wedge \tau}\left(L_R^2(s)|Y_s^{\varepsilon}|^2+|\alpha(s,X^{\e}(s),0,q_s)|^2+\|\beta(s,X^{\e}(s),0,q_s)\|^2\right)ds,\\
<& \infty.
\end{split}
\end{equation*}
By the proof of Theorem \ref{solutionX} we have that $\mathbb{E}\sup_{t\in [0,T]}|X_{t}^{\varepsilon}|<\infty$, and then
\begin{equation}
\label{eqn:Yassmprob}
\begin{split}
\mathbb{P}\left(\sup_{t\in [0,T]}|Y_{t}^{\varepsilon}|
> M \right)&=\mathbb{P}\left(\sup_{t\in [0,T]}|Y_{t}^{\varepsilon}|
> M, T\leq \tau\right)+\mathbb{P}\left(\sup_{t\in [0,T]}|Y_{t}^{\varepsilon}|
> M, T>\tau\right)\\
&\leq \mathbb{P}\left(\sup_{t\in [0,T]}|Y_{t\wedge \tau}^{\varepsilon}|
> M\right)+\mathbb{P}(T>\tau)\\
&\leq \frac{\mathbb{E}\left(\sup_{t\in [0,T]}|Y_{t\wedge \tau}^{\varepsilon}|^2\right)}{M^2}+\mathbb{P}\left(\sup_{t\in [0,T]}|X_{t}^{\varepsilon}|> R\right)\\
&\xrightarrow[]{M\to\infty ~\mbox{and then} ~R\rightarrow \infty}0.
\end{split}
\end{equation}
We further define another stopping time $\bar{\tau}$ as
\begin{equation}
\label{eqn:bartau}
\bar{\tau}=\tau \wedge \inf\{s: |Y_s^{\varepsilon}|>M\}.
\end{equation}
Then by It\^o's formula and the Gr{\"o}nwall's lemma, we have
\begin{align*}
&\mathbb{E}\sup_{t\in [0,T\wedge \bar \tau]}|Y_{t}^{\varepsilon}-Y_{t}|^2\\
\leq & C\mathbb{E}\int_0^{t \wedge \bar \tau} \left|\alpha(s,X_s^{\varepsilon},Y_s^{\varepsilon},q_s)-\alpha(s,X_s,Y_s,q_s)\right|^2ds\\
&+ C\mathbb{E}\int_0^{t \wedge \bar \tau} \left\|\beta(s,X_s^{\varepsilon},Y_s^{\varepsilon},q_s)-\beta(s,X_s,Y_s,q_s)\right\|^2ds\\
\leq & C \mathbb{E} \int_0^{t \wedge \bar \tau} \left(L_R^2(s)|Y_s^{\varepsilon}-Y_s|^2+|\alpha(s,X_s^{\varepsilon},Y_s,q_s)-\alpha(s,X_s,Y_s,q_s)|^2\right)ds\\
&+C \mathbb{E} \int_0^{t \wedge \bar \tau} \|\beta(s,X_s^{\varepsilon},Y_s,q_s)-\beta(s, X_s,Y_s,q_s)\|^2 ds\\
\leq & C \mathbb{E} \int_0^{t \wedge \bar \tau} |\alpha(s,X_s^{\varepsilon},Y_s,q_s)-\alpha(s,X_s,Y_s,q_s)|^2 ds\\
& + C \mathbb{E} \int_0^{t \wedge \bar \tau}  \|\beta(s,X_s^{\varepsilon},Y_s,q_s)-\beta(s,X_s, Y_s,q_s)\|^2 ds.
\end{align*}
Then by the continuity of functions $\alpha$ and $\beta$ enforced in Condition \ref{Y}, as well as the the convergence result in Theorem \ref{Xpathasymp}, we have 
$$\mathbb{E}\sup_{t\in [0,T\wedge \bar \tau]}|Y_{t}^{\varepsilon}-Y_{t}|^2\to  0,\quad\quad\quad\mbox{as}\;\; \varepsilon\to0.$$
By equation \eqref{eqn:Yassmprob} and similar to its derivation, we can obtain that for any $\eta>0$,
$$\lim_{\varepsilon\to0}\mathbb{P}\left(\sup_{t\in [0,T]}|Y_t^{\varepsilon}-Y_t|>\eta\right)= 0.$$
\end{proof}

\section{One-dimensional SVI  system with H\"older continuous coefficients}
\label{sec:holder}

In this section, we consider the following one-dimensional SVI  system with H\"older continuous coefficients:
\begin{equation}
\label{eqn:holder}
\left\{\begin{array}{lll}
 X_t\in & x_0+\int_0^tb(s,X_s)ds+\int_0^t\sigma_1(s,X_s)dW_s+\int_0^t\sigma_2(s, X_s) dB_s\\
 &-\int_0^t\partial\psi_1(X_s)ds,\\
       \\
Y_t\in & y_0+\int_0^t \alpha(s,X_s,Y_s, q_s)ds+\int_0^t\beta(s,X_s,Y_s, q_s)dB_s\\
&-\int_0^t\partial\psi_2(Y_s)ds,\\
             \end{array}\right.
\end{equation}
where $b, ~\sigma_1, ~\sigma_2$ are measurable functions mapping from $\mR^+\times\mR$ to $\mR$, $\alpha$ and $\beta$ are measurable functions mapping from $\mR^+\times\mR\times\mR\times \mathbb{U}$ to $\mR$, $W$ and $B$ are two independent standard one-dimensional Brownian motions on a complete filtered probability space $(\Omega, \mathcal {F}, \mathcal {F}_t, \mathbb{P})$.

\begin{condition}\label{1dX} For the $X$ process in the SVI  system \eqref{eqn:holder}, we impose the following conditions: Assume that $b(t,x), ~\sigma_1(t,x), ~\sigma_2(t,x)$ are continuous in $(t,x)$, and
\begin{align}
&\left(b(t, x)-b(t, x')\right)(x-x')\leq 0\nonumber\\
&\left(b(t, x)-b(t, x')\right)^2\leq l_0(t) (x-x')^{1+2\alpha}, \quad&\text{for some } \alpha\in[0,1/2],\nonumber\\
&\left(\sigma_i(t, x)-\sigma_i(t, x')\right)^2\leq l_i(t)(x-x')^{1+2\alpha}, \quad &i=1, 2,\nonumber\\
& \psi_1\geq\psi_1(0)=0, \quad & 0\in\mathrm{Int}(D_1),\nonumber
\end{align}
where $l_i(t)$ for $i=0,1,2$ are functions of $t$ only and satisfy $l_i(\cdot)\in L^1([0,T])$.
\end{condition}

\begin{condition}\label{1dY} For the $Y$ process in the SVI  system \eqref{eqn:holder}, we impose the following conditions:
\begin{itemize}
\item $\lambda_1\leq q_t\leq \lambda_2$,
\item $\alpha, ~\beta$ are continuous in $(t,x,y,q)$ satisfying
$$(y-y')\big(\alpha(t,x,y,q)-\alpha(t,x,y',q)\big)\leq 0,$$
and for $\gamma\in[0,1/2]$
\begin{equation*}\begin{split}
&|\alpha(t,x,y,q)-\alpha(t,x',y',q)|^2\vee|\beta(t,x,y)-\beta(t,x',y')|^2\\
\leq&c(t)(|x-x'|^{1+2\gamma}+|y-y'|^{1+2\gamma}),
\end{split}\end{equation*}
where $c(t)$ is locally integrable for any $t\geq0$,.
\item $0\in\mathrm{Int}(D_2)$, $\psi_2\geq\psi_2(0)\equiv0$.
\end{itemize}
\end{condition}

\subsection{Well-posedness}
\label{sec:holder_wellposedness}

First of all we solve the well-posedness problem under the above conditions. An estimate for the solution process is given in the following proposition.

\begin{proposition}\label{1dproperty1}
Suppose $(X,\phi^{(1)})$ is a solution of the $X$ process in the SVI  system \eqref{eqn:holder}, under Condition \ref{1dX}, one has
$$\mathbb{E}\|X\|_T^2+\mathbb{E}\int_0^T\psi_1(X_s)ds\leq C(1+|x_0|^2),$$
and then
$$\mathbb{E}|\phi^{(1)}|_T^0\leq C(1+|x_0|^2).$$
\end{proposition}

\begin{proof}
Note that by Condition \ref{1dX},
\begin{equation}\label{bsigma}
\begin{split}
&|b(t,x)|^2\leq 2|b(t,x)-b(t,0)|^2+|b(t,0)|^2\leq l_0(t)|x|^2+|b(t,0)|^2.\\
&|\sigma_i(t,x)|^2\leq l_i(t)|x|^2+|\sigma_i(t,0)|^2, ~~~i=1,2.
\end{split}
\end{equation}
Then applying It\^o's formula and by Remark \ref{remark2.1}, we have
\begin{align*}
|X_t|^2=&|x_0|^2+2\int_0^tX_sb(s,X_s)ds+2\int_0^tX_s\sigma_1(s,X_s)dW_s+2\int_0^tX_s\sigma_2(s,X_s)dB_s\\
&+\sum_{i=1}^2\int_0^t|\sigma_i(s,X_s)|^2ds-2\int_0^tX_sd\phi^{(1)}_s\\
\leq&|x_0|^2+\int_0^t\big(1+l_0(s)+l_1(s)+l_2(s)\big)|X_s|^2ds\\
&+\int_0^t\big(|b(s,0)|^2+|\sigma_1(s,0)|^2+|\sigma_2(s,0)|^2\big)ds\\
&+2\int_0^t(X_s,\sigma_1(s,X_s))dW_s+2\int_0^t(X_s,\sigma_2(s,X_s))dB_s-2\int_0^t\psi_1(X_s)ds.
\end{align*}
By using the BDG inequality and the H\"older's inequality,
\begin{align*}
&\mathbb{E}\sup_{t\leq T}\left|2\int_0^tX_s\sigma_1(s,X_s)dW_s+2\int_0^tX_s\sigma_2(s,X_s)dB_s\right|\\
\leq&C\mathbb{E}\left(\int_0^T|X_s|^2|\sigma_1(s,X_s)|^2ds\right)^{1/2}+C\mathbb{E}\left(\int_0^T|X_s|^2|\sigma_2(s,X_s)|^2ds\right)^{1/2}\\
\leq&C\mathbb{E}\|X\|_T^2+C\mathbb{E}\int_0^T\big(1+l_1(s)+l_2(s)\big)|X_s|^2ds+C\mathbb{E}\int_0^T\big(|\sigma_1(s,0)|^2+|\sigma_2(s,0)|^2\big)ds.
\end{align*}
Therefore, Gr{\"o}nwall's lemma yields that
\begin{equation*}
\mathbb{E}\|X\|_T^2+\mathbb{E}\int_0^T\psi_1(X_s)ds\leq C(1+|x_0|^2).
\end{equation*}
Moreover, by using this estimate and Remark \ref{remark2.1}, we also have
\begin{equation*}
\mathbb{E}|\phi^{(1)}|_T^0\leq C(1+|x_0|^2).
\end{equation*}
\end{proof}

The well-posedness of the $X$ process in the SVI  system \eqref{eqn:holder} is established in the following proposition.
\begin{proposition}\label{well1d}
Under Condition \ref{1dX}, there is a unique strong solution of the $X$ process in the SVI  system \eqref{eqn:holder}.
\end{proposition}

\begin{proof}
We apply a regularization approximation method here. Define the Moreau-Yosida regularization of $\psi_1$ as
\be\label{MY}
\psi_1^n(x):=\inf \left\{\frac n{2}|x'-x|^2+\psi_1(x'); x'\in\mR \right\}, ~~~n\geq1, ~~\forall x\in\mR.
\ee
Then $\psi_1^n$ is a $\mathcal{C}^1$-convex function, and its gradient $\nabla\psi_1^n$ is monotone and Lipschitz with Lipschitz constant $n$ which is due to the reason that $\nabla\psi_1$ has no gradient.
Moreover, according to \cite{asiminoaei1997approximation}, $\nabla\psi_1^n$ has the following properties 
\begin{align}
&(x-x')(\nabla\psi_1^n(x)-\nabla\psi_1^m(x'))\geq-\left(\frac1n+\frac1m\right)\nabla\psi_1^n(x)\nabla\psi_1^m(x'),~~~\forall x, x'\in\mathbb R, \label{psinmproperty1}\\
&\nabla\psi_1^n(x)\in\partial\psi_1(J_nx),\quad \psi_1(J_nx)\leq \psi_1^n(x)\leq\psi_1(x),\label{psinmproperty2}\\
&\psi_1^n(x)=\psi_1^n(J_nx)+\frac1{2n}|\nabla\psi_1^n(x)|^2, \label{psinmproperty3}
\end{align}
where $J_nx:=x-\frac1n\nabla\psi_1^n(x)$.

It is known that the following stochastic differential equation has a unique strong solution 
\begin{equation}
\label{eqn:holder_perturbed_system}
 dX_t^n= b(t,X_t^n)dt+\sigma_1(t,X_t^n)dW_t+\sigma_2(t, X_t^n) dB_t-\nabla\psi_1^n(X_t^n)dt,~~~
 X_0^n=x_0\in \bar{D}_1,
 \end{equation}
where $\nabla\psi_1^n$ is the gradient of $\psi_1^n$.

Moreover, with arguments similar to those in Proposition \ref{1dproperty1},
\begin{align*}
\mathbb{E}\|X^n\|_T^4\leq&C\mathbb{E}|x_0|^4+C\mathbb{E}\Big(\int_0^T\big(1+l_0(s)+l_1(s)+l_2(s)\big)|X^n_s|^2 ds\Big)^2\\
&+C\mathbb{E}\Big(\int_0^T\big(|b(s,0)|^2+|\sigma_1(s,0)|^2+|\sigma_2(s,0)|^2\big)ds\Big)^2\\
&+\mathbb{E}\Big(\int_0^T(X^n_s,\sigma_1(s,X^n_s))dW_s+\int_0^T(X^n_s,\sigma_2(s,X^n_s))dB_s\Big)^2\\
\leq&C\mathbb{E}|x_0|^4+\frac12\mathbb{E}\|X^n\|_T^4+C\mathbb{E}\left(\int_0^T\big(1+l_0(s)+l_1(s)+l_2(s)\big)|X^n_s|^2 ds\right)^2\\
&+C\mathbb{E}\big(\int_0^T\big(|b(s,0)|^2+|\sigma_1(s,0)|^2+|\sigma_2(s,0)|^2\big)ds\big)^2\\
\leq&C\mathbb{E}|x_0|^4+\frac12\mathbb{E}\|X^n\|_T^4+C_T\mathbb{E}\left(\int_0^T\big(1+l_0(s)+l_1(s)+l_2(s)\big)^2|X^n_s|^4ds\right)\\
&+C\mathbb{E}\big(\int_0^T\big(|b(s,0)|^2+|\sigma_1(s,0)|^2+|\sigma_2(s,0)|^2\big)ds\big)^2,
\end{align*}
where in the last inequality we used the Cauchy-Schwarz inequality in the integral form.
Then Gr\"onwall's lemma yields
\begin{equation}
\label{4moment}
\sup_n\mathbb{E}\|X^n\|_T^4\leq C(1+\mathbb{E}|x_0|^4),
\end{equation}
and by the dynamic \eqref{eqn:holder_perturbed_system} we further have
\begin{equation}\label{psin}
\sup_n\mathbb{E}\left(\int_0^T|\nabla\psi_1^n(X_s^n)|ds\right)^2<\infty.
\end{equation}

Note that by It\^o's formula, the fact that $\nabla\psi_1^n$ is Lipschitz with Lipschitz constant $n$, and equation \eqref{psinmproperty3}, we have
\begin{align*}
&|\psi_1^n(X^n_t)|^2\\
=&|\psi_1^n(x_0)|^2+2\int_0^t\psi_1^n(X^n_s)\nabla\psi_1^n(X^n_s)b(s,X^n_s)ds-2\int_0^t\psi_1^n(X^n_s)|\nabla\psi_1^n(X^n_s)|^2ds\\
&+\sum_{i=1}^2\int_0^t|\nabla\psi_1^n(X^n_s)|^2|\sigma_i(s,X^n_s)|^2ds+n\sum_{i=1}^2\int_0^t\psi_1^n(X^n_s)|\sigma_i(s,X^n_s)|^2ds\\
&+2\int_0^t\psi_1^n(X^n_s)\nabla\psi_1^n(X^n_s)\sigma_1(s,X^n_s)dW_s+2\int_0^t\psi_1^n(X^n_s)\nabla\psi_2^n(X^n_s)\sigma_2(s,X^n_s)dB_s\\
\leq&|\psi_1^n(x_0)|^2+2n\int_0^t\psi_1^n(X^n_s)|X^n_s b(s,X^n_s)|ds-2\int_0^t\psi_1^n(X^n_s)|\nabla\psi_1^n(X^n_s)|^2ds\\
&+3n\sum_{i=1}^2\int_0^t\psi_1^n(X^n_s)|\sigma_i(s,X^n_s)|^2ds+2\int_0^t\psi_1^n(X^n_s)\nabla\psi_1^n(X^n_s)\sigma_1(s,X^n_s)dW_s\\
&+2\int_0^t\psi_1^n(X^n_s)\nabla\psi_1^n(X^n_s)\sigma_2(s,X^n_s)dB_s.
\end{align*}
By the BDG's inequality, Condition \ref{1dX}, equation \eqref{psinmproperty3}, and the Young's inequality for products, we obtain
\begin{equation*}\begin{split}
&\mathbb{E}\sup_{t\leq T}\left|2\int_0^t\psi_1^n(X^n_s)\nabla\psi_1^n(X^n_s)\sigma_1(s,X^n_s)dW_s+2\int_0^t\psi_1^n(X^n_s)\nabla\psi_1^n(X^n_s)\sigma_2(s,X^n_s)dB_s\right|\\
\leq&C\mathbb{E}\Big(\int_0^T\left|\psi_1^n(X^n_s)\nabla\psi_1^n(X^n_s)\sigma_1(s,X^n_s)\right|^2ds\Big)^{1/2}\\
&+C\mathbb{E}\Big(\int_0^T|\psi_1^n(X^n_s)\nabla\psi_1^n(X^n_s)\sigma_2(s,X^n_s)|^2ds\Big)^{1/2}\\
\leq&\frac12\mathbb{E}\sup_{t\leq T}|\psi_1^n(X^n_t)|^2+Cn\mathbb{E}\int_0^T|\psi_1^n(X^n_s)|\cdot \big(|\sigma_1(s,X^n_s)|^2+|\sigma_2(s,X^n_s)|^2\big)ds.
\end{split}\end{equation*}
By the fact that $$|\psi_1^n(X^n_s)|\leq |\nabla\psi_1^n(X^n_s)|\cdot |X^n_s|$$ since $\psi_1^n$ is a convex function, and by the Young's inequality for products, we have 
\begin{align*}
&\frac12\mathbb{E}\sup_{t\leq T}|\psi_1^n(X^n_t)|^2+2\mathbb{E}\int_0^T\psi_1^n(X^n_s)|\nabla\psi_1^n(X^n_s)|^2ds\\
\leq&C\mathbb{E}|\psi_1^n(x_0)|^2+Cn\mathbb{E}\int_0^T|\psi_1^n(X^n_s)|\bigg(|X^n_s||b(s,X^n_s)|+|\sigma_1(s,X^n_s)|^2+|\sigma_2(s,X^n_s)|^2\bigg)ds\\
\leq&C\mathbb{E}|\psi_1^n(x_0)|^2+Cn\mathbb{E}\int_0^T|\psi_1^n(X^n_s)|^{1/3}|\nabla\psi_1^n(X^n_s)|^{2/3}|X^n_s|^{2/3}\\
&\hspace{120pt minus 1fil}\times\bigg(|X^n_s||b(s,X^n_s)|+|\sigma_1(s,X^n_s)|^2+|\sigma_2(s,X^n_s)|^2\bigg)ds\hfilneg\\
\leq&C\mathbb{E}|\psi_1^n(x_0)|^2+\mathbb{E}\int_0^T|\psi_1^n(X^n_s)||\nabla\psi_1^n(X^n_s)|^2ds\\
&+Cn^{3/2}\mathbb{E}\int_0^T|X^n_s|\left(|X^n_s|^2+\sum_{i=0}^2l_i(s)|X^n_s|^2+|b(s,0)|^2+|\sigma_1(s,0)|^2+|\sigma_2(s,0)|^2\right)ds\\
\leq&C\mathbb{E}|\psi_1^n(x_0)|^2+\mathbb{E}\int_0^T|\psi_1^n(X^n_s)||\nabla\psi_1^n(X^n_s)|^2ds\\
&+Cn^{3/2}\mathbb{E}\int_0^T(1+|X^n_s|^4)\left(1+\sum_{i=0}^2l_i(s)+|b(s,0)|^2+|\sigma_1(s,0)|^2+|\sigma_2(s,0)|^2\right)ds,
\end{align*}
which together with equation \eqref{4moment} yields that
\begin{equation}
\mathbb{E}\sup_{t\leq T}|\psi_1^n(X^n_t)|^2\leq Cn^{3/2}.
\end{equation}
By equation \eqref{psinmproperty3}, we further have
\begin{equation}\label{psin4}
\mathbb{E}\sup_{t\leq T}|\nabla\psi_1^n(X^n_t)|^4\leq 4n^2\mathbb{E}\sup_{t\leq T}|\psi_1^n(X^n_t)|^2\leq Cn^{7/2}.
\end{equation}

Now take any $\delta \in (0,1)$, any $h>0$, and set
$$g_{\delta,h}(x)=\int_0^x\int_0^y f_{\delta,h}(\gamma)d\gamma dy$$
where $f_{\delta,h}\geq 0$ and vanishes outside $[h\delta,h]$, and
$$f_{\delta,h}(x)\leq \frac{2}{x\ln \delta^{-1}}, \quad \int f_{\delta,h}(x) dx=1.$$
Then we have
\begin{equation}
\label{eqn:delta_h_inequality}
|x|\leq g_{\delta,h}(|x|)+h,
\end{equation}
and
\begin{equation}
\label{eqn:delta_h_inequality_2}
0\leq {g}_{\delta,h}'\leq 1, \quad {g}_{\delta,h}''(|x|)\leq \frac{2}{x\ln \delta^{-1}}\mathbbm{1}_{(|x|\in [h\delta,h])}.
\end{equation}

By applying equation \eqref{eqn:delta_h_inequality} and then It\^o's formula, 
\begin{align*}
|X^{m}_t-X^n_t|\leq&g_{\delta,h}(|X^{m}_t-X^n_t|)+h\\
\leq&\int_0^tg'_{\delta,h}(|X^{m}_s-X^n_s|)\frac{X^{m}_s-X^n_s}{|X^{m}_s-X^n_s|}\big[b(s,X^{m}_s)-b(s,X^n_s)\big]ds\\
&+\frac12\sum_{i=1}^2\int_0^tg''_{\delta,h}(|X^{m}_s-X^n_s|)\big[\sigma_i(s,X^{m}_s)-\sigma_i(s,X^n_s)\big]^2ds\\
&+\int_0^tg'_{\delta,h}(|X^{m}_s-X^n_s|)\frac{X^{m}_s-X^n_s}{|X^{m}_s-X^n_s|}\big[\sigma_1(s,X^{m}_s)-\sigma_1(s,X^n_s)\big]dW_s\\
&+\int_0^tg'_{\delta,h}(|X^{m}_s-X^n_s|)\frac{X^{m}_s-X^n_s}{|X^{m}_s-X^n_s|}\big[\sigma_2(s,X^{m}_s)-\sigma_2(s,X^n_s)\big]dB_s\\
&-\int_0^tg'_{\delta,h}(|X^{m}_s-X^n_s|)\frac{X^{m}_s-X^n_s}{|X^{m}_s-X^n_s|}\big[\nabla\psi_1^m(X^{m}_s)-\nabla\psi_1^n(X^n_s)\big]ds+h.
\end{align*}
Then by Condition \ref{1dX}, equation \eqref{psinmproperty1}, and equation  \eqref{eqn:delta_h_inequality_2},
$$|X^{m}_t-X^n_t|\leq I(t)+M(t)+J(t)+h,$$
where 
\begin{align*}
I(t):=&\frac{1}{\ln \delta^{-1}}\sum_{i=1}^2\int_0^tl_i(s)|X^{m}_s-X^n_s|^{2\alpha}\mathbbm{1}_{\{|X^{m}_s-X^n_s|\in[h\delta,h]\}}ds,\\
M(t):=&\int_0^tg'_{\delta,h}(|X^{m}_s-X^n_s|)\frac{X^{m}_s-X^n_s}{|X^{m}_s-X^n_s|}\big[\sigma_1(s,X^{m}_s)-\sigma_1(s,X^n_s)\big]dW_s\\
&+\int_0^tg'_{\delta,h}(|X^{m}_s-X^n_s|)\frac{X^{m}_s-X^n_s}{|X^{m}_s-X^n_s|}\big[\sigma_2(s,X^{m}_s)-\sigma_2(s,X^n_s)\big]dB_s,\\
J(t):=&\int_0^tg'_{\delta,h}(|X^{m}_s-X^n_s|)\left(\frac1n+\frac1m\right)|X^{m}_s-X^n_s|^{-1}\nabla\psi_1^m(X^{m}_s)\nabla\psi_1^n(X^n_s)ds.
\end{align*}
Clearly we have 
$$\mathbb{E}\sup_{t\leq T}|I(t)|\leq \frac{2h^{2\alpha}}{\ln \delta^{-1}}\sum_{i=1}^2\int_0^tl_i(s)\leq C\frac{h^{2\alpha}}{\ln \delta^{-1}},$$
and 
\begin{equation*}\begin{split}
\mathbb{E}\sup_{t\leq T}|M(t)|\leq& C\sum_{i=1}^2\mathbb{E}\left(\int_0^Tl_i(s)|X^{m}_s-X^n_s|^{1+2\alpha}ds\right)^{1/2}\\
\leq&C\sum_{i=1}^2\mathbb{E}\int_0^Tl_i(s)|X^{m}_s-X^n_s|^{2\alpha}ds+\frac12\mathbb{E}\sup_{t\leq T}|X^{m}_t-X^n_t|.
\end{split}\end{equation*}
By equations \eqref{psin}, \eqref{psin4}, and \eqref{eqn:delta_h_inequality_2},
\begin{align*}
\mathbb{E}\sup_{t\leq T}|J(t)|\leq& \frac{1}{h\delta}\mathbb{E}\int_0^T \left(\frac1n+\frac1{m}\right)\nabla\psi_1^{m}(X^{m}_s)\nabla\psi_1^n(X^n_s)ds\\
\leq&\frac{1}{h\delta}\Big[\frac1n\big(\mathbb{E}\sup_{t\leq T}|\nabla\psi_1^n(X^n_t)|^2\big)^{1/2}\Big(\mathbb{E}\big(\int_0^T|\nabla\psi_1^{m}(X^{m}_t)|dt\big)^{2}\Big)^{1/2}\\
&+\frac1m\big(\mathbb{E}\sup_{t\leq T}|\nabla\psi_1^m(X^{m}_t)|^2\big)^{1/2}\Big(\mathbb{E}\big(\int_0^T|\nabla\psi_1^n(X^n_t)|dt\big)^{2}\Big)^{1/2}\Big]\\
\leq&C\frac1{h\delta}(n^{-1/8}+m^{-1/8}).
\end{align*}
Summing up these estimates, by Gr\"onwall's lemma, we have
\begin{equation*}\begin{split}
\mathbb{E}\sup_{t\leq T}|X^{m}_t-X^n_t|\leq& C\frac1{h\delta}(n^{-1/8}+m^{-1/8})+C\frac{h^{2\alpha}}{\ln \delta^{-1}}+h.
\end{split}\end{equation*}
Considering $\alpha\in[0,1/2]$, we further have
\begin{equation*}
\mathbb{E}\sup_{t\leq T}|X^{m}_t-X^n_t| \leq C(h\delta)^{-2\alpha}(n^{-\alpha/4}+m^{-\alpha/4})+C\frac{h^{2\alpha}}{\ln \delta^{-1}}+h.
\end{equation*}
Taking $\delta=\frac12$ and $h=\min\{m,n\}^{-16}$ yields
\begin{equation*}
\mathbb{E}\sup_{t\leq T}|X^{m}_t-X^n_t|\leq C\min\{m,n\}^{-\frac{\alpha}8}\to0, \quad\mbox{as}\quad n\to\infty.
\end{equation*}
Moreover, by setting $$\phi^{(1),n}_t:=\int_0^t\nabla\psi_1^n(X^n_s)ds,$$ we have
\begin{equation*}
\mathbb{E}\sup_{t\leq T}|\phi^{(1),m}_t-\phi^{(1),n}_t|\to0, \quad\mbox{as}\quad n\to\infty.
\end{equation*}
Hence $(X^n, \phi^{(1),n})$ is Cauchy in the complete metric space $$L^1(\Omega;\mathcal{C}([0,T];\mR))\times L^1(\Omega;\mathcal{C}([0,T];\mR))$$ and thus there exists $(X,\phi^{(1)})$ in the space satisfying that
\begin{equation}
\label{eqn:Xn_cvg}
\mathbb{E}\sup_{t\leq T}|X^{n}_t-X_t|\to0 \quad\text{and}\quad \mathbb{E}\sup_{t\leq T}|\phi^{(1),n}_t-\phi^{(1)}_t|\to0,\quad\quad\mbox{as}\quad n\to\infty.
\end{equation}

Now it remains to prove that $(X,\phi^{(1)})$ is a solution. Since by equations \eqref{psin} we have
$$
\sup_n\mathbb{E}\|\phi^{(1),n}\|_T<\infty,
$$
it then yields
$$
\mathbb{E}\|\phi^{(1)}\|_T<\infty.
$$
Recall that $\psi_1^n$ is convex and that
$$
\psi_1(J_nx)\leq\psi_1^n(x)\leq \psi_1(x)
$$
given in equation \eqref{psinmproperty2}, for any $\varrho\in \mathcal{C}([0,T];\mR)$ and any $t\in[0,T]$,
\begin{equation*}\begin{split}
\int_0^t(\varrho_s-X^n_s)d\phi^{(1),n}_s=&\int_0^t(\varrho_s-X^n_s)\nabla\psi_1^n(X^n_s)ds\\
\leq& \int_0^t\psi_1^n(\varrho_s)ds-\int_0^t\psi_1^n(X^n_s)ds\\
\leq&\int_0^t\psi_1^n(\varrho_s)ds-\int_0^t\psi_1(J_nX^n_s)ds.
\end{split}\end{equation*}
By equation \eqref{eqn:Xn_cvg} and the fact that monotone increasing sequence of random variables that converge in probability implies convergence almost surely, sending $n\to\infty$ gives
$$
\int_0^t(\varrho_s-X_s)d\phi^{(1)}_s\leq \int_0^t\psi_1(\varrho_s)ds-\int_0^t\psi_1(X_s)ds.
$$
Hence $(X,\phi^{(1)})$ is a solution.

\end{proof}

With analogous arguments, we can obtain that there exists a unique strong solution for the $Y$ process in the SVI  system \eqref{eqn:holder} and the proof is omitted.

\subsection{Asymptotic analysis}
\label{sec:holder_asymptotic}
In this section, we perform asymptotic analysis on the perturbed one-dimensional SVI  system \eqref{eqn:holder_perturbed_stystem} with H\"older continuous coefficients described in Condition \ref{1dX_asymptotic}, regarding its limiting system \eqref{eqn:holder} satisfying Condition \ref{1dX}.

The perturbed version of the one-dimensional SVI  system \eqref{eqn:holder} with a small positive parameter $\epsilon$ is given by
\begin{equation}
\label{eqn:holder_perturbed_stystem}
\left\{\begin{array}{lll}
 X_t^{\varepsilon}\in & x_0+\int_0^tb^{\varepsilon}(s,X_s^{\varepsilon},\varepsilon)ds+\int_0^t\sigma_1^{\varepsilon}(s,X_s^{\varepsilon},\varepsilon)dW_s+\int_0^t\sigma_2^{\varepsilon}(s, X_s^{\varepsilon}, \varepsilon) dB_s\\
 &-\int_0^t\partial\psi_1(X^{\varepsilon}_s)ds,\\
       \\
Y_t^{\varepsilon}\in & y_0+\int_0^t \alpha(s,X_s^{\varepsilon},Y_s^{\varepsilon}, q_s)ds+\int_0^t\beta(s,X_s^{\varepsilon},Y_s^{\varepsilon}, q_s)dB_s\\
&-\int_0^t\partial\psi_2(Y^{\varepsilon}_s)ds,
\end{array}\right.
\end{equation}
where 
\begin{equation}
\label{eqn:barX_functions}
\begin{split}
\lim_{\varepsilon \rightarrow 0} b^{\varepsilon}(t, x, \varepsilon)=b(t, x), \quad \lim_{\varepsilon \rightarrow 0} \sigma_i^{\varepsilon}(t, x, \varepsilon)=\sigma_i(t, x), \quad i=1, 2.
\end{split}
\end{equation}

\begin{condition}\label{1dX_asymptotic} Assume that $b^{\varepsilon}(t,x,\varepsilon), ~\sigma_1^{\varepsilon}(t,x,\varepsilon), ~\sigma_2^{\varepsilon}(t,x,\varepsilon)$ are continuous in $(t,x)$, uniformly in $\varepsilon$, and
\begin{align}
&\left(b^{\varepsilon}(t, x, \varepsilon)-b^{\varepsilon}(t, x', \varepsilon)\right)(x-x')\leq 0,\nonumber\\
&\left(b^{\varepsilon}(t, x, \varepsilon)-b^{\varepsilon}(t, x', \varepsilon)\right)^2\leq l_0(t) (x-x')^{1+2\alpha}, \quad&\text{for some } \alpha\in[0,1/2],\nonumber\\
&\left(\sigma_i^{\varepsilon}(t, x, \varepsilon)-\sigma_i^{\varepsilon}(t, x', \varepsilon)\right)^2\leq l_i(t)(x-x')^{1+2\alpha}, \quad &i=1, 2,\nonumber\\
& \psi_1\geq\psi_1(0)=0, \quad & 0\in\mathrm{Int}(D_1),\nonumber
\end{align}
where $l_i(t)$ for $i=0,1,2$ are functions of $t$ only and satisfy $l_i(\cdot)\in L^1([0,T])$.
\end{condition}


With arguments similar to those of Proposition \ref{1dproperty1}, one can obtain the following proposition.
\begin{proposition}
\label{prop:sup_2ndmoment}
Under Conditions \ref{1dX} and \ref{1dX_asymptotic}, one has
\begin{equation}
\mathbb{E} \sup_{t\in [0,T]}|X_t|^2<\infty\quad\text{and}\quad \sup_{\varepsilon}\mathbb{E} \sup_{t\in [0,T]}|X_t^{\varepsilon}|^2<\infty.
\end{equation}
\end{proposition}

In the following, we give the convergence result regarding the $X_t^{\varepsilon}$ process as $\varepsilon$ goes  to $0$.
\begin{proposition}
Under Conditions \ref{1dX} and \ref{1dX_asymptotic}, as $\varepsilon \rightarrow 0$, we have
\begin{equation}
\begin{split}
\mathbb{E}\sup_{t\in [0,T]}|X_t^{\varepsilon}-X_t|\rightarrow 0.
\end{split}
\end{equation}
\end{proposition}

\begin{proof}
By equation \eqref{eqn:delta_h_inequality} and It\^o's formula, one has
\begin{align}
|X_t^{\varepsilon}-X_t|\leq & g_{\delta,h}(|X_t^{\varepsilon}-X_t|)+h\nonumber\\
=&\int_0^t {g}_{\delta,h}'(|X_s^{\varepsilon}-X_s|)\frac{X_s^{\varepsilon}-X_s}{|X_s^{\varepsilon}-X_s|}\big[b^{\varepsilon}(s, X_s^{\varepsilon}, \varepsilon)-b(s,  X_s)\big]ds\nonumber\\
&+\frac12\int_0^t {g}_{\delta,h}''(|X_s^{\varepsilon}-X_s|)\sum_{i=1}^2\big[\sigma_i^{\varepsilon}(s, X_s^{\varepsilon}, \varepsilon)-\sigma_i(s, X_s)\big]^2ds\nonumber\\
&+\int_0^t {g}_{\delta,h}'(|X_s^{\varepsilon}-X_s|)\frac{X_s^{\varepsilon}-X_s}{|X_s^{\varepsilon}-X_s|}\big[\sigma_1^{\varepsilon}(s, X_s^{\varepsilon}, \varepsilon)-\sigma_1(s, X_s)\big]dW_s\label{eqn:Xcvg_expansion}\\
&+\int_0^t {g}_{\delta,h}'(|X_s^{\varepsilon}-X_s|)\frac{X_s^{\varepsilon}-X_t}{|X_s^{\varepsilon}-X_s|}\big[\sigma_2^{\varepsilon}(s, X_s^{\varepsilon}, \varepsilon)-\sigma_2(s, X_s)\big]dB_s+h.\nonumber
\end{align}
Note that by Condition \ref{1dX_asymptotic} and that $g'_{\delta,h}\in[0,1]$,
\begin{equation}
\begin{split}
&\int_0^t {g}_{\delta,h}'(|X_s^{\varepsilon}-X_s|)\frac{X_s^{\varepsilon}-X_s}{|X_s^{\varepsilon}-X_s|}\big[b^{\varepsilon}(s, X_s^{\varepsilon}, \varepsilon)-b(s, X_s)\big]ds\\
= & \int_0^t {g}_{\delta,h}'(|X_s^{\varepsilon}-X_s|)\frac{X_s^{\varepsilon}-X_s}{|X_s^{\varepsilon}-X_s|}\big[b^{\varepsilon}(s, X_s^{\varepsilon}, \varepsilon)-b^{\varepsilon}(s, X_s, \varepsilon)+b^{\varepsilon}(s, X_s, \varepsilon)
-b(s, X_s)\big]ds\\
\leq & \int_0^t {g}_{\delta,h}'(|X_s^{\varepsilon}-X_s|)\frac{X_s^{\varepsilon}-X_s}{|X_s^{\varepsilon}-X_s|}|b^{\varepsilon}(s, X_s, \varepsilon)-b(s, X_s)|ds\\
\leq & \int_0^t |b^{\varepsilon}(s, X_s, \varepsilon)
-b(s, X_s)|ds,
\end{split}\end{equation}
and by Condition \ref{1dX} and Proposition \ref{1dproperty1} one has
\begin{align*}
&\sup_{\varepsilon}\mathbb{E}\left(\int_0^T |b^{\varepsilon}(s, X_s^{\varepsilon}, \varepsilon)-b(s,X_s)|ds\right)^2\\
\leq &C \sup_{\varepsilon}\mathbb{E}\int_0^T |b^{\varepsilon}(s, X_s^{\varepsilon}, \varepsilon)|^2ds+C\sup_{\varepsilon}\mathbb{E}\int_0^T|b(s,X_s)|^2ds\\
\leq& C \sup_{\varepsilon}\mathbb{E}\int_0^T\big(l_0(s)|X_s^{\varepsilon}|^{1+2\alpha}+|b^{\varepsilon}(s,0,\varepsilon)|^2+|b(s,0)|^2\big)ds\\
<&\infty.
\end{align*}
Hence, by equation \eqref{eqn:barX_functions}, as $\varepsilon\rightarrow 0$,
\begin{align*}
\mathbb{E}\int_0^t {g}_{\delta,h}'(|X_s^{\varepsilon}-X_s|)\frac{X_s^{\varepsilon}-X_s}{|X_s^{\varepsilon}-X_s|}\big[b^{\varepsilon}(s, X_s^{\varepsilon}, \varepsilon)-b(s,  X_s)\big]ds\rightarrow 0.
\end{align*}
Similarly, by Propositions \ref{1dproperty1} and \ref{prop:sup_2ndmoment} as well as the regularity conditions for $\sigma_i$ and $\sigma_i$ respectively, one has
\begin{align*}
&\sup_{\varepsilon}\mathbb{E}\int_0^T\big[\sigma_i^{\varepsilon}(s,  X_s, \varepsilon)-\sigma_i(s, X_s)\big]^2ds\\
\leq &C \sup_{\varepsilon}\mathbb{E}\int_0^T\big[\sigma_i^{\varepsilon}(s, X_s, \varepsilon)-\sigma_i^{\varepsilon}(s, 0, \varepsilon)]^2ds+C\sup_{\varepsilon}\mathbb{E}\int_0^T[\sigma_i^{\varepsilon}(s, 0, \varepsilon)-\sigma_i(s, 0)\big]^2ds\\
&+C\sup_{\varepsilon}\mathbb{E}\int_0^T\big[\sigma_i(s, X_s)-\sigma_i(s, 0)\big]^2ds\\
<&\infty.
\end{align*}

Then, by equation \eqref{eqn:barX_functions}, as $\varepsilon\to0$,
\begin{equation}
\mathbb{E}\int_0^T\big[\sigma_i^{\varepsilon}(s, X_s, \varepsilon)-\sigma_i(s, X_s)\big]^2ds\to0,~~~~i=1,2,
\end{equation}
and for sufficiently small $\varepsilon$ satisfying that
$$
\mathbb{E}\int_0^T\big[\sigma_i^{\varepsilon}(s, X_s, \varepsilon)-\sigma_i(s, X_s)\big]^2ds<\delta h^{1+2\alpha},
$$ 
by equation \eqref{eqn:delta_h_inequality_2} one has
\begin{align*}
&\mathbb{E}\int_0^t {g}_{\delta,h}''(|X_s^{\varepsilon}-X_s|)\big[\sigma_i^{\varepsilon}(s, X_s^{\varepsilon}, \varepsilon)-\sigma_i(s,  X_s)\big]^2ds\\
\leq & C\mathbb{E}\int_0^t \frac{l_i(s)}{\ln \delta^{-1}|X_s^{\varepsilon}-X_s|}|X_s^{\varepsilon}-X_s|^{1+2\alpha}\mathbbm{1}_{\{|X_s^{\varepsilon}-X_s|\in[h\delta,h]\}}ds\\
&+C\mathbb{E}\int_0^t \frac{1}{\ln \delta^{-1}|X_s^{\varepsilon}-X_s|}\big[\sigma_i^{\varepsilon}(s X_s, \varepsilon)-\sigma_i(s, X_s)\big]^2\mathbbm{1}_{\{|X_s^{\varepsilon}-X_s|\in[h\delta,h]\}}ds\\
\leq & C\mathbb{E}\int_0^t \frac{l_i(s)}{\ln \delta^{-1}}h^{2\alpha}ds+C\mathbb{E}\int_0^t \frac{1}{\delta h \ln \delta^{-1}}\big[\sigma_i^{\varepsilon}(s, X_s, \varepsilon)-\sigma_i(s,  X_s)\big]^2ds\\
\leq & \frac{C h^{2\alpha}}{\ln \delta^{-1}}+\frac{C}{\delta h \ln \delta^{-1}}\mathbb{E}\int_0^t\big[\sigma_i^{\varepsilon}(s, X_s, \varepsilon)-\sigma_i(s, X_s)\big]^2ds\\
\leq & \frac{C h^{2\alpha}}{\ln \delta^{-1}}.
\end{align*}
Plugging the above results in equation \eqref{eqn:Xcvg_expansion}, one gets
$$\mathbb{E}|X_t^{\varepsilon}-X_t|\leq \frac{C h^{2\alpha}}{\ln \delta^{-1}}+h.$$
Taking supremum and then expectation of equation \eqref{eqn:Xcvg_expansion}, we obtain
\begin{align*}
&\mathbb{E}\sup_{t\in [0,T]}|X_t^{\varepsilon}-X_t|\\
\leq & \frac{C h^{2\alpha}}{\ln \delta^{-1}}+h+C\sum_{i=1}^2\mathbb{E}\left(\int_0^T \big[\sigma_i^{\varepsilon}(t, X_s^{\varepsilon}, \varepsilon)-\sigma_i(s, X_s)\big]^2ds\right)^{1/2}\\
\leq & \frac{C h^{2\alpha}}{\ln \delta^{-1}}+h+C\sum_{i=1}^2\mathbb{E}\left(\int_0^T l_i(s) |X_s^{\varepsilon}- X_s|^{1+2\alpha}ds\right)^{1/2}\\
&+C\sum_{i=1}^2\mathbb{E}\left(\int_0^T \big[\sigma_i^{\varepsilon}(t, X_s, \varepsilon)-\sigma_i(s, X_s)\big]^2ds\right)^{1/2}\\
\leq & \frac{C h^{2\alpha}}{\ln \delta^{-1}}+h+\frac{1}{2}\mathbb{E}\sup_{t\in [0,T]}|X_t^{\varepsilon}-X_t|+C\sum_{i=1}^2\mathbb{E}\int_0^T l_i(s) |X_s^{\varepsilon}- X_s|^{2\alpha}ds\\
&+C\sum_{i=1}^2\mathbb{E}\left(\int_0^T \big[\sigma_i^{\varepsilon}(t, X_s, \varepsilon)-\sigma_i(s, X_s)\big]^2ds\right)^{1/2},
\end{align*}
where we used the H\"older's inequality in the last equality to reduce the order of $|X_s^{\varepsilon}- X_s|$ on the right hand side. By reorganizing the terms and noticing that $2\alpha<1$, Gr\"onwall's lemma yields
$$\mathbb{E}\sup_{t\in [0,T]}|X_t^{\varepsilon}-X_t|\leq \frac{C h^{2\alpha}}{\ln \delta^{-1}}+h.$$
Taking $h=\delta$, we have
$$\mathbb{E}\sup_{t\in [0,T]}|X_t^{\varepsilon}-X_t|\leq  \frac{C \delta^{2\alpha}}{\ln \delta^{-1}}+\delta.$$
For the reason that $\delta$ can take any small values, we conclude the proof as desired.
\end{proof}

With arguments analogous to Theorem \ref{thm:Ycvg}, we also have the following convergence result for the $Y$-system whose proof is omitted.
 \begin{theorem} Under Conditions \ref{1dY} and \ref{1dX_asymptotic}, one has
$$
\mathbb E\sup_{t\leq T}|Y^{\varepsilon}_t-Y_t|^2\to 0, ~~~\mbox{as}~~\varepsilon\to0.
$$

\end{theorem}



\bibliography{references}
\end{document}